\documentclass[10pt,oneside,amscd,amssymb,epsf]{amsart}

\usepackage{amsthm}
\usepackage{graphicx}
\usepackage{pb-diagram, Bdiagx}
\usepackage{psfrag}

\newcommand{\TryPackage}[3]{\IfFileExists{#1.sty}{\usepackage{#1}
#2}{#3}}
\TryPackage{mathrsfs}{\renewcommand{\mathcal}{\mathscr }}{%
        \TryPackage{eucal}{}{}}

\begin{document}
\title[Calderon Projector for the Hessian ...]{Calderon Projector for
the
Hessian of the perturbed Chern-Simons function on a 3-manifold with
boundary}

\author{Benjamin Himpel}
\address{Department of Mathematics, Indiana University, Bloomington, IN
47405}
\email{bhimpel@indiana.edu}
\urladdr{http://mypage.iu.edu/$\sim$bhimpel}

\author{Paul Kirk}
\address{Department of Mathematics, Indiana University, Bloomington, IN
47405}
\email{pkirk@indiana.edu}
\urladdr{http://mypage.iu.edu/$\sim$pkirk}

\author{Matthias Lesch}
\address{Universit\"at zu K\"oln, Mathematisches Institut, Weyertal
86--90,
D--50931 K\"oln}
\email{lesch@mi.uni-koeln.de}
\urladdr{http://www.mi.uni-koeln.de/$\sim$lesch}

\newcommand{\cormark}{\marginpar{$\Leftarrow$}}
\date{Feb. 18, 2003}
\begin{abstract} {The existence and continuity for the Calder\'on
projector of the perturbed odd signature  operator on a 3-manifold is
established. As an application we give a new proof of a result of Taubes
relating the mod 2 spectral flow of a family of operators on a  homology
3-sphere  with the difference in local intersection numbers of the
character varieties  coming from a Heegard decomposition. }
\end{abstract}

\thanks{The second named author gratefully acknowledges the support of
the
National Science Foundation under   grant no. DMS-9971020.}

\maketitle

\renewcommand{\tilde}{\widetilde}
\newcommand{\Cald}{Calder\'on projector}
\newcommand{\cV}{{\mathcal V}}

\newcommand{\const}{\mbox{const}}
\newcommand{\lto}{\longrightarrow}
\newcommand{\al}{\alpha}
\newcommand{\be}{\beta}
\newcommand{\ga}{\gamma}
\newcommand{\de}{\delta}
\newcommand{\ep}{\epsilon}
       \renewcommand{\th}{\theta}
\newcommand{\la}{\lambda}
\newcommand{\La}{\Lambda}
\newcommand{\om}{\omega}
\newcommand{\si}{\sigma}
\newcommand{\Ga}{\Gamma}
\newcommand{\Om}{\Omega}
\newcommand{\Si}{\Sigma}
\newcommand{\etab}{\tilde{\eta}}
\newcommand{\tsig}{\tilde{\sigma}}
\newcommand{\hal}{\hat{\alpha}}

\newcommand{\grad}{\mbox{grad}}
\newcommand{\gra}{\mbox{grad}^\sim}

\newcommand\tr{\mathop {\rm tr}}

\newcommand{\tlam}{\widetilde \lambda}

\newcommand{\zz}{{\mathbb Z}}
\newcommand{\rr}{{\mathbb R}}
\newcommand{\cc}{{\mathbb C}}
\newcommand{\qq}{{\mathbb Q}}
\newcommand{\hh}{{\mathbb H}}

\newcommand{\cA}{{\mathcal A}}
\newcommand{\cU}{{\mathcal U}}
\newcommand{\cE}{{\mathcal E}}
\newcommand{\cD}{{\mathcal D}}
\newcommand{\cC}{{\mathcal C}}
\newcommand{\cL}{{\mathcal L}}
\newcommand{\cB}{{\mathcal B}}
\newcommand{\cF}{{\mathcal F}}
\newcommand{\cG}{{\mathcal G}}
\newcommand{\cM}{{\mathcal M}}
\newcommand{\cH}{{\mathcal H}}
\newcommand{\cT}{{\mathcal T}}
\newcommand{\cO}{{\mathcal O}}
\newcommand{\cP}{{\mathcal P}}
\newcommand{\hA}{\widehat{A}}
\newcommand{\hB}{\widehat{B}}

\newcommand{\tB}{{\widetilde{\mathcal B}}}
\newcommand{\tM}{{\widetilde{\mathcal M}}}
\newcommand{\tC}{{\widetilde{C}}}

\newcommand{\im}{\operatorname{im}}
\newcommand{\Map}{\operatorname{Map}}
\newcommand{\hol}{\operatorname{{\it hol}}}
\newcommand{\proj}{\operatorname{Proj}}
\newcommand{\Mas}{\operatorname{Mas}}
\newcommand{\Spec}{\operatorname{Spec}}
\newcommand{\Hom}{\operatorname{Hom}}
\newcommand{\Hol}{\operatorname{Hol}}
\newcommand{\crit}{\operatorname{Crit}}
\newcommand{\hess}{\operatorname{Hess}}
\newcommand{\ind}{\operatorname{ind}}
\newcommand{\image}{\operatorname{image}}
\newcommand{\ad}{\operatorname{ad}}
\newcommand{\End}{\operatorname{End}}
\newcommand{\Sign}{\operatorname{Sign}}
\newcommand{\SF}{\operatorname{SF}}
\newcommand{\Id}{\operatorname{Id}}
\def\eqref#1{(\ref{#1})}

\newcommand{\del}{\partial}

\newcommand{\chern}{\operatorname{ch}}
\renewcommand{\eqref}[1]{\textnormal{(\ref{#1})}}

       \newcommand{\opX}{{\bar{X}}}

   \theoremstyle{definition}

\newtheorem{defn}{Definition}[section]
\theoremstyle{plain}
\newtheorem{lemma}[defn]{Lemma}
\newtheorem{thm}[defn]{Theorem}
\newtheorem{prop}[defn]{Proposition}
       \newtheorem{cor}[defn]{Corollary}
     \theoremstyle{definition}
\newtheorem{remark}[defn]{Remark}

\numberwithin{equation}{section}

\section{Introduction}

The purpose of this article is to establish the existence and continuity
of the Calder\'on projectors for families of perturbed signature
operators
on a 3-manifold with boundary.  We explain briefly the motivation.

The Morse theory  of the Chern-Simons function on the space of $SU(n)$
connections over a 3-manifold leads to the construction of topological
invariants, notably Taubes' construction of Casson's invariant
\cite{taubes}, Floer's instanton homology \cite{floer}, and more
recently
$SU(3)$ and $SU(n)$ extensions of Casson's invariant \cite{herald-boden,
BHK-int, CLM, herald2}.  These invariants are defined using the Morse
index of the Hessian of the Chern-Simons function at a critical point.
To
make this rigorous Taubes showed (in the $SU(2)$ setting) how to
understand this ill-defined Morse index in terms of the spectral flow of
the odd signature operator coupled to a path of connections. Most
importantly for the present article, he also showed how to construct
perturbations of the Chern-Simons function in order to make it suitably
non-degenerate.

Taubes' ideas form the blueprint for   generalizations of Casson's
invariant. In particular the definition of the generalized Casson
invariants involve the spectral flow of 1-parameter families of
perturbed
odd signature operators.  These perturbed odd signature operators (which
we describe in detail below) are compact perturbations of Dirac
operators,
but they are not differential operators.

With the construction of a Morse-theoretic topological invariant
established, the problem arises of computing the invariant to discover
what information it contains.  Casson's surgery formula and Floer's
exact
triangle are examples of computational formulas which make it possible
to
compute their invariants. In general one wishes to have surgery formulas
for the generalized Casson invariants.  This leads to the problem of
computing the spectral flow of the Hessian of the Chern-Simons function between
  critical points in terms of cut--and--paste  data for the underlying
3-manifold.

There is a good cut--and--paste theory for the spectral flow of Dirac
operators which make it possible to carry out  calculations.
Nicolaescu's article \cite{nico} gives a splitting formula for the
spectral flow of a family of Dirac operators. This formula equates the
spectral flow to an infinite--dimensional Maslov index corresponding to
the
associated paths of Calder\'on projectors coming from two parts in a
decomposition of the manifold. (Roughly speaking, the \Cald\  for a
Dirac
operator $D$ is the projection in
$L^2(\partial X)$ onto the subspace of boundary values of solutions to
the equation $D\phi=0$.)
    Nicolaescu's result, combined with some techniques introduced in
\cite{daniel-kirk}, sufficed to compute $SU(3)$ Casson invariants for
certain Seifert fibered homology 3-spheres in
\cite{BHKK,CLM}, where one could avoid using perturbations (at least
near
the reducible critical points). However, the presence of singularities in the flat moduli space 
complicates calculations for more general
homology spheres, requiring the use of perturbations. To apply the powerful
methods of \cite{nico} one needs an extension of the results concerning
the existence and continuity of the Calder\'on projectors for families
of
operators more general than Dirac operators; in fact general enough to
include the perturbed odd signature operator. This is the task
accomplished in this article. Our main result is the following, which
combines Theorems
\ref{continuityCald} and \ref{Poissonop-a}.

\bigskip

\noindent{\bf Theorem.} {\em Consider the class $\cD$ of Dirac operators
on $E\to X$ for a manifold with boundary $X$ which are in cylindrical
form
near
$\partial X$,  and the class $\mathcal{V}$ of auxiliary operators
satisfying the pseudolocality property
\textnormal{\eqref{pseudolocality}}.  Then for $D\in \cD$ and $V\in
\mathcal{V}$ the Calder\'on projector of $D+V$ is well defined.

We say that a family $D_t\in\cD$ varies continuously if all coefficients
of $D_t$ (in any local chart) vary continuously. Moreover, we equip
$\mathcal{V}$ with the usual norm topology of $\cB(L^2(X,E))$. Then for
each
$n\in\mathbb{Z}_+$ the map
\[ \bigl\{(D,V)\in\cD\times\mathcal{V}\,\big|\, \dim Z_0(  D+
V)=n\bigr\}
      \longrightarrow \cB(L^2(\Sigma,E)), \quad (D,V)\mapsto C(  D+
V)
\] taking a pair $(D,V)$ to the  Calder\'on projector $C$  for the operator
$D+V$ is continuous. }

\smallskip

  In the statement, $Z_0(D+V)$   is a vector space isomorphic to $K\oplus
K$, where
$K$ denotes the space of sections $\phi$ which vanish on a collar of the
boundary and which satisfy  $(D+V)\phi=0$.  For Dirac operators the
unique continuation property implies that
$Z_0=0$. A novel feature of our approach is that we do not require
the unique continuation property to hold.

\bigskip
 
We are chiefly concerned with the application of this result to the
perturbed odd signature operator arising from the Morse theory of the
Chern-Simons function, as this will be  used crucially in a forthcoming
article \cite{BHK-pqr}. (In fact the completion of that article was held
up by  the need for the main result of the present article.)
Unfortunately
there is no suitably precise exposition of the various constructions
which
arise in this context and so we include  the details necessary  to
set up
the problem  precisely and also for the convenience of the reader.  One
of our
purposes  is to explicitly provide justification for comments like ``in
the non-generic
case one must perturb, but the arguments carry through,''    in the
context
of cut--and--paste techniques for spectral flow.

It is our hope that in addition to providing the proof of the theorem
stated
above this article complements  the  exposition \cite{Herald1, herald2}.
Thus in Section
\ref{set-up} we  explain the approach used in the Morse-theoretic study
of
the Chern-Simons function, with an emphasis on   how the perturbations
are constructed and on the effect  of perturbing on the Hessian of the
Chern-Simons function. We also include proofs of   some  (known) results
concerning the pseudolocality of the perturbation (Proposition
\ref{lem2.1}) and an estimate concerning the  size of the perturbation
(Lemma \ref{estimate}) that we require  to apply the result stated above
to the context of the  Hessian of the perturbed Chern-Simons function.

Theorem \ref{thm5.3} is our desired application. In this statement,
$D_{A,f}$ denotes the
perturbed twisted odd signature operator coupled to a connection $A$
and a perturbation $f$. It is
essentially the Hessian of the perturbed Chern-Simons function
$cs+f:\cA\to \mathbb{R}$.

\bigskip

\noindent{\bf Theorem \ref{thm5.3}.} {\em Let $X$ be a compact
3-manifold
with boundary and $\cA$ the space of   connections in a principal
$SU(n)$
bundle over $X$. Then there exists a neighborhood
$U\subset \cA$ of the space of flat connections, and an
$\epsilon>0$  so that if
$B_{\cP}(0,\epsilon)$ denotes the
$\epsilon$ ball around $0$ in the space of perturbations $\cP$,  the map
$U\times B_{\cP}(0,\epsilon)\to \cB(L^2(\Sigma,E))$ sending a pair
$(A,f)$ to the Calder\'on projector $P_{A,f}$ of $D_{A,f}$   is
continuous.}

\bigskip

We finish the article   by giving a new proof in Section
\ref{application}, using our methods, of the main step in Taubes'
argument
that his invariant
equals Casson's invariant.

\bigskip

\noindent{\bf Theorem} (Taubes){\bf .} {\em The mod 2 Spectral flow of
the
perturbed odd signature operator between two generic perturbed flat
connections on a homology 3-sphere equals zero if and only if the
corresponding  local intersection numbers of the $SU(2)$ character
varieties of the Heegard handlebodies in the $SU(2)$ character variety
of
the Heegard surface are equal. }

\bigskip

  The strategy of our proof is to use an adiabatic stretching of a collar
of the Heegard surface  and symplectic reduction to reduce the
infinite--dimensional Maslov index (which equals the sign Taubes assigns
to
a critical point) to a finite--dimensional Maslov index in  the character
variety of a Heegard surface  (which determines Casson's sign). Having
established in  Theorem
\ref{thm5.3} that the Calder\'on projectors exist and vary  nicely, the
first step in the argument uses this and Nicolaescu's theorem
\cite{nico} to identify spectral flow with  the infinite--dimensional
Maslov index of the  Calder\'on projectors. The rest of the argument is
a study of the effect of stretching and a corresponding symplectic
reduction,  as analyzed in
\cite{daniel-kirk}.

The authors thank  H. Boden, G. Daskalopoulos, C. Herald, T. Mrowka,
and K.
Wojciechowski for helpful discussions.

\section{basic set-up}\label{set-up}
\subsection{Connections on a 3-manifold}

We begin by recalling the  spaces and operators of interest.

Let $X$ be a compact 3-manifold.  Fix a principal $SU(n)$ bundle  $P$
over
$X$. We associate two vector bundles to $P$: the $\cc^n$ vector bundle
$E\to X$ associated   to the standard representation
$SU(n)\to GL(\cc^n)$, and the Lie algebra bundle
$\ad P\to X$ associated to the adjoint representation $SU(n)\to
GL(su(n))$
of
$SU(n)$ on its Lie algebra $su(n)$.  Notice that the bundle $P$ is
trivializable since
$\pi_1(SU(n))$ and $\pi_2(SU(n))$  are trivial.

We let $\cA$ denote the space of smooth connections on $P$. A connection
$A\in \cA$ defines  a covariant derivative
$$d_A:C^\infty(X;E)\to \Om^1(X;E).$$ It is  convenient to view $\ad
P$
as a subbundle of  $\End(E)$. Then the covariant derivative on $\ad P$
associated to $A$  (which we also denote by $d_A:C^\infty(X,\ad P)\to
\Om^1(X;\ad P)$)  satisfies
$$(d_A(\phi))(e)=d_A(\phi(e))-\phi(d_A(e))\  \text{ for } \
\phi\in C^\infty(X;\ad P), \ e\in C^\infty(X;E).$$ The covariant
derivatives    extend  to $p$-forms $d_A:\Om^p(X; E)\to
\Om^{p+1}(X; E)$  and  $d_A:\Om^p(X; \ad P)\to
\Om^{p+1}(X; \ad P)$ by requiring the graded Leibnitz rule to hold. The
operator $d^2_A:\Omega^p(X;E)\to\Omega^{p+2}(X;E)$ has order zero and
is
given by multiplication by the   curvature form  $F(A)\in
\Omega^2(X;ad P)$.

     The space
$\cA$ is an affine space modeled  on  the space of Lie algebra valued
1-forms
$\Om^1(X;\ad P)$: Fixing a connection $A_0$ yields an identification
\begin{equation}\label{eq2.1}\cA\xrightarrow\cong \Om^1(X;\ad P),
\ A\mapsto A-A_0.\end{equation} It is sometimes convenient to  fix a
trivialization $P\cong X\times SU(n)$; then we can take $A_0$ to be the
trivial (product) connection which we denote by $\Theta$. The space
$\cA$
is acted on by the   group of (smooth) gauge transformations
$\cG=C^\infty(X;P\times_{Ad} SU(n))$ in the usual way: $d_{gA}(e)=
g(d_A(g^{-1}(e)))$.

We fix a Riemannian metric on $X$. If the boundary of $X$ is non-empty
we
will always assume that a collar of the boundary $\del X$ in $X$ is
isometric to $[0,1)\times \del X$.

The $L^2$ inner product on
$\Om^i(X;\ad P)$ is then defined in terms of the Hodge $*$ operator on
$X$ by the formula
\begin{equation}\label{eq2.2}\langle b,c\rangle_{L^2}=-\int_X tr(b\wedge
*c).\end{equation}

For analytical reasons one must at times complete the spaces $\cA$ and
$\cG$ in appropriate Sobolev norms.
     Thus we fix some $s\ge 1$ and take
$\cA_s$ to be the
$L^2_s$ completion of $\cA\cong \Om^1(X;\ad P)$, and $\cG_{s+1}$ to be
the
closure of
$\cG$ in the
$L^2_{s+1}$ completion of $C^\infty(\End(E))$.  Here $L^2_s$ denotes
     the completion in the Sobolev norm corresponding to sections with
$s$
derivatives in $L^2$, extended in the standard way to all
$s\in \rr$. Then
$\cG_{s+1}$ acts smoothly on $\cA_{s}$ (see \cite{FU}). The restriction
$s\ge 1$ guarantees that the gauge transformations are continuous.

\subsection{The Chern-Simons function} The Chern-Simons function
$cs:\cA\to
\rr$ is defined in terms of a trivialization $P\cong X\times SU(n)$ as
follows. To a connection
$A\in \cA$ one can assign the $su(n)$-valued 1-form
$$a:=A-\Theta\in
\Om^1(X;\ad P)\cong\Om^1(X)\otimes su(n).$$ One then defines
$$cs(A)= \frac{1}{8\pi^2}\int_X \tr(da \wedge a +\tfrac{2}{3} a\wedge
a\wedge a).$$ In this expression (as in Equation \eqref{eq2.2}) the
wedge
product means the wedging of the differential forms and matrix
multiplication of the
$su(n)$  coefficients. Note that the definition (and value) of the
Chern-Simons function depends on the choice of trivialization of
$P$. If $X$ is {\it closed}, then the reduction of $cs(A)$ modulo
$\zz$ is independent of the choice of trivialization, and is gauge
invariant. In fact, $cs(g\cdot A)= cs(A) + \text{degree}(g) $, where the
degree is defined using the fact that
$H^3(SU(n);\zz)=\zz$.  Thus when $X$ is closed, one considers the
Chern-Simons function as a circle valued function on the quotient
$cs:\cA/\cG\to \rr/\zz=S^1$.

In \cite{taubes}, Taubes gave an interpretation of Casson's invariant of
homology 3-spheres by studying the Morse theory of the Chern-Simons
function. This was generalized in much subsequent work, including the work of
Floer
\cite{floer} and Herald-Boden \cite{herald-boden}.
Thus we look next   at the differential topology of
$\cA$,
$\cG$,  and $cs$.

The tangent space of $\cA$ at a connection $A$ is canonically identified
with $\Om^1(X;\ad P)$.  It is not hard to identify the tangent space of
$\cG$  at the identity with $\Om^0(X;\ad P)$. Fixing a connection $A\in
\cA$ yields  the {\em action map}
$$g_A:\cG\to \cA, \  g\mapsto g\cdot A$$ The differential of this map is
just the map $d_A:\Om^0(X;\ad P)\to \Om^1(X;\ad P)$. Precisely, the
diagram
\begin{equation}\label{eq2.3}
\begin{diagram}\dgmag{900}\dgsquash
\node{T_1\cG}\arrow{e,t}{dg_A}\arrow{s,l}{\cong}
          \node{T_A\cA}\arrow{s,l}{\cong}\\
\node{\Om^0(X;\ad P)}\arrow{e,t}{d_A}\node{\Om^1(X;\ad P)}
\end{diagram}\end{equation} commutes. Thus the tangent space to the
quotient $\cA/\cG$ at $[A]$,
$T_{[A]}\cA/\cG$ can be identified with the cokernel of $d_A$ (at least
at
smooth points of $\cA/\cG$: making this precise requires looking at
Sobolev completions and dealing with quotient singularities but we will
not emphasize this point here).

\subsection{The Hessian of the Chern-Simons function}\label{sec2.3} {\em
In this section
we assume that $X$ is a closed 3-manifold}.  This is because the
differential operators which arise from the Morse theory of the
Chern-Simons function  have a simple  description on a closed manifold
which is adequate for our purposes, even though we are ultimately
interested in cutting and pasting questions. This section could be
written
in the context of manifolds with boundary, but in that case the proper
way
to study the Chern-Simons function is as a section of a determinant
bundle
as in \cite{RSW} and \cite{Herald1}. (The material of this subsection is
well known. Standard references include \cite{taubes,
donaldson-kronheimer, FU}.   We include it to set up notation and to
prepare the reader for
the section on perturbations.)

The differential
$dcs_A:T_A\cA\to
\rr$ of the Chern-Simons function is computed (see \cite{taubes}) {\em
for
$X$ closed} as:
\[ dcs_A(b)=\frac{1}{4\pi^2}\int_X \tr(F(A)\wedge b),\] with $F(A)$
the curvature 2-form of $A$. (Recall that
in
a  trivialization,
     $F(A) =da+a\wedge a$.)
     From this one concludes that the critical points of $cs:\cA\to\rr$
are
those
$A$ such that $F(A)=0$, i.e. the {\em flat connections}.

The Hessian of $cs:\cA\to\rr$, $\hess_A(cs):T_A\cA\times T_A\cA\to
\rr$  is defined for any connection  $A\in\cA$ using the affine
structure
of $\cA$. It can be computed explicitly when
$X$ is closed (\cite{taubes}):
\[
\hess(cs)_A(b,c)=-\frac{1}{4\pi^2}\langle *d_Ab,c\rangle.
\]

The $L^2$ inner product can be viewed as a Riemannian metric on the
infinite--dimensional manifold $\cA$ (or $\cA_s$). Identifying
vectors and covectors via 
this metric  allows us to define the gradient vector field of
$cs$
     by the rule $\langle\grad(cs)_A,b\rangle_{L^2}=dcs_A(b)$. Similarly
the
linearization $H_A:T_A\cA\to T_A\cA$ of $\grad(cs)$ is obtained
from the Hessian by the rule
$\langle H_A(b),c\rangle_{L^2}=\hess(cs)_A(b,c)$.  From the above
formulas
one concludes that $$\grad(cs)_A= -\tfrac{1}{4\pi^2}*F(A)$$ and
$$H_A(b)= -\tfrac{1}{4\pi^2}*d_A(b).$$

The self--adjoint operator $H_A$ is not elliptic because of the gauge
invariance of
$cs$. For example, at a flat connection $A$, $d_Ad_A=F(A)=0$ , and so
$H_A(d_A(b))=0$ for all $b\in \Om^0(X;\ad P)$. In other words the kernel
of $H_A$ contains the infinite--dimensional subspace of tangent vectors
to
the
$\cG$ orbit
     through $A$. This causes problems if one wishes to consider the
Morse
theory of $cs$.

As explained in \cite{taubes} (and is typical in gauge theory), the
correct
way to resolve this lack of ellipticity is to consider $\grad(cs)$
as a
vector field   on the quotient
$\cA/\cG$, and thus $H_A$ as a family of elliptic self--adjoint
operators on
$T_{[A]}\cA/\cG$.    Using the slice theorem  and Diagram
\eqref{eq2.3} one identifies
$T_{[A]}\cA/\cG$ (one has to take care at the singular points) as the
orthogonal complement of the image of
$d_A:\Om^0(X;\ad P)\to \Om^1(X;\ad P)$, in other words as the kernel of
the adjoint $d_A^*:\Om^1(X;\ad P)\to \Om^0(X;\ad P)$.  The Bianchi
identity implies that
$d_A^*(\grad(cs)(A))=-\tfrac{1}{4\pi^2}d_A^*(*F(A))=0$. If $A$ is a
critical point of $cs$ (i.e.~a flat connection), then
$$d_A^*(H_A(b))=\tfrac{1}{4\pi^2}*d_A**d_A(b)=
\tfrac{1}{4\pi^2}*d_Ad_A(b)=\tfrac{1}{4\pi^2}*F(A)(b)=0.$$
But $H_A(b)$
is
not necessarily in the kernel of $d_A^*$ when $A$ is not a critical
point.
This motivates setting
$$\tilde{H}_A :\ker d_A^*\to \ker d_A^*, \ \tilde{H}_A(b)=
-\frac{1}{4\pi^2}\proj_{\ker d_A^*} (*d_A(b))$$ where $\proj$ denotes the $L^2$
orthogonal
projection. The operator $\tilde{H}_A$ agrees with the restriction of
$H_A$ to $\ker d_A^*$ if the connection $A$ is flat.

Then $\tilde{H}_A$ is  Fredholm  (i.e. its extension to
$L^2_s(\Omega^1(X;\ad P))\cap \ker d_A^*\to L^2_{s-1}(\Omega^1(X;\ad
P))\cap \ker d_A^*$ is Fredholm). Moreover, it is self--adjoint with
compact resolvent with respect to the $L^2$ inner product restricted to
$\ker d_A^*$. The spaces $\ker d_A^*$ form a smooth vector bundle over
each orbit-type stratum of
$\cA/\cG$. In particular, Taubes' idea was to use  the spectral flow
   of
a
family of operators $\tilde{H}_{A_t}$ as a substitute for the difference
in Morse index of critical points $A_0$ and $A_1$ of the Chern-Simons
function, and thereby to define an Euler characteristic for
$\cA/\cG$.

(The {\em spectral flow} $\SF(H_t)$ of a 1-parameter family of
self--adjoint operators $H_t$ is
the integer defined to be the   number of eigenvalues  that cross zero, counted with sign.
Some conventions need to be set if the
operators in question have kernel at the endpoints of the path, and one
needs to
establish that the concept is well defined  for appropriately
continuous  paths of
unbounded operators. A careful construction appropriate for our
purposes can can be
found in
\cite{booss-lesch-phillips}. The basic properties and conventions are
explained in
\cite{kirk-lesch}.)

Taubes observed that the spectral flow of a family $\tilde{H}_{A_t}$
equals the spectral flow of a slightly better behaved family of
operators.
Consider, for any connection $A$, the twisted de Rham sequence
\begin{equation}\label{eq2.4} 0\to \Om^0(X;\ad
P)\xrightarrow{d_A}\Om^1(X;\ad P)\xrightarrow{d_A}\Om^2(X;\ad
P)\xrightarrow{d_A}\Om^3(X;\ad P)\to 0.
\end{equation} The sequence \eqref{eq2.4} is an  elliptic  complex when
$A$ is flat.   It is not a complex at a non-flat connection, but can be
made into a complex if one substitutes
$d'_A:=-4\pi^2 *\tilde{H}_A\circ\proj_{\ker
d_A^*}$ for
$d_A:\Omega^1(X;\ad P)\to \Omega^2(X;\ad P)$.   When $A$ is flat,
$d'_A=d_A$.

Whether or not
$A$ is flat, the   {\it odd signature operator}  (obtained essentially
by folding up the sequence \eqref{eq2.4})
     \begin{equation}
\label{foldedupop}\begin{split} &D_A:\Om^0(X;\ad P)\oplus \Om^1(X;\ad
P)\to \Om^0(X;\ad P)\oplus
\Om^1(X;\ad P)\\ &D_A(b,c)=(d_A^* c, *d_A c+d_A b)
\end{split}\end{equation} is a self--adjoint Dirac operator in the
sense of   e.g. \cite{booss-w}
  for any
connection $A$.

One also has a  self--adjoint operator $D'_A(b,c)=(d_A^* c, *d'_A c+
d_Ab)$. For $A$ flat, $D_A=D'_A$. For a general irreducible connection
$A$,
$D_A-D'_A$ is a   relatively compact operator. 
In fact, letting $\Delta_A^{-1}:\Om^0(X;\ad P)\to \Om^0(X;\ad P)$ denote the inverse of $d_A^*d_A$  (which is invertible for $A$ irreducible) and observing   that
$\proj_{\ker
d_A^*}x= x-d_A\Delta_A^{-1}d_A^*x$, we see that 
$$(D_A-D'_A)(b,c)= \big( *F(A)\Delta_A^{-1} d_A^*- d_A\Delta_A^{-1}*F(A)+d_A\Delta_A^{-1}*F(A)d_A\Delta_A^{-1}d_A^*\big)(c).$$
 (see \cite{taubes} for details).

  By decomposing
$\Om^0(X;\ad P)\oplus
\Om^1(X;\ad P)$ as $\Om^0(X;\ad P)\oplus  \text{Image }d_A\oplus\ker
d_A^*
$, one sees that   the spectrum of
$D'_A$ is the disjoint union of the spectrum of $\tilde{H}_A$ and a
symmetric (with respect to $\lambda\to-\lambda$) spectrum.
A simple
argument then shows that given a path $A_t$ of (say irreducible)
connections joining two critical points of $cs$,   the spectral flow of
the family
$\tilde{H}_{A_t}$ (which plays the role of the difference in  Morse
indices of $cs$ at $A_0$ and
$A_1$) equals the spectral flow of the family $D_{A_t}$. What is gained
is that the spectral flow of $D_{A_t}$ makes sense whether or not $A_t$
is a path of irreducible connections, i.e.~even when the subspaces
$\ker d_A^*\subset \Omega^1(X; \ad P)$ do not vary continuously.

\subsection{Perturbations} We come now to the construction of
perturbations and their effect on
the Hessian of the Chern-Simons  function.
   The critical point set of
$cs$ may be complicated and it is often necessary to perturb the
Chern-Simons function to make it as nice as possible.

     Making precise what ``nice" means can be quite difficult in the case
of
$SU(n)$ for $n>2$, or when the homology of $X$ is complicated. In the
case
of $SU(2)$ and $X$ a homology sphere  Taubes replaced $cs$ by a function
$cs+h$ whose restriction to the top stratum of
$\cA/\cG$ (i.e. the irreducible connections; see Section
\ref{application})   is Morse in a suitable sense, and then used this to
define a topological invariant in the spirit of the Poincar\'e-Hopf
theorem,  which he showed equals Casson's invariant. In
\cite{herald-boden} a similar analysis is carried out for $SU(3)$ and
$X$
a homology sphere leading to a definition of an $SU(3)$ Casson
invariant;
in that case in addition to being Morse on the top stratum of
$\cA/\cG$,
$cs+h$ must be suitably non-degenerate in higher codimensional strata.
In
more recent work \cite{herald2} Herald shows how with a suitable choice
of
perturbation
$h$ one can define versions of Casson's invariant for rational homology
spheres and $SU(n)$, $n=2,3,4,5$.

There are slightly different approaches to constructing a suitable
family
of perturbations $h$, but they are similar and our method will apply to
each approach.  We mention that the basic requirements that $h$ should
satisfy are:
\begin{enumerate}
\item The perturbation $h$ should be a smooth function on the Hilbert
manifold $\cA_1/\cG_2$. 
\item The operator 
defined by the Hessian of $cs +h$  at $A$ should be a
compact perturbation of
$\tilde{H}_A$.
\item The family of admissible $h$ should be large enough so that one
can
prove various general position results about $cs+h$ on the strata and
their normal bundles of $\cA/\cG$.
\end{enumerate}

The following types of perturbations,   due to Taubes and Floer,
accomplish these goals, as explained in
\cite{herald2}.

First, fix a collection of smooth embeddings $\gamma_i:D^2\times S^1\to
X,\ i=1,\cdots,N$. Taubes \cite{taubes} requires that the embeddings
have
disjoint images, whereas   Floer \cite{floer} and Herald
\cite{herald2}  assume that   $\gamma_i(x,1)=\gamma_1(x,1)$ for all
$i$ and $x\in D^2$ and that the derivatives of the $\gamma_i$ in the
$S^1$
direction  at points in $D^1\times \{1\}$ are independent of $i$. (In
fact, nothing is lost by assuming that there is an interval $I$ around
$1\in S^1$ so that $\gamma_i(x,u)$ is independent of $i$ for $u\in I$.)
The results of this article work equally well in   each setting. The
difference between them is that it is easier to show that the third
condition above holds for   perturbations of the type Floer considers
since  the base point issues are easier to control. We will assume, to
keep the notation under control, that we are using the approach of
Floer,
i.e. the common basepoint approach.
    We do not address the third condition in this article.  In the
language
of
\cite{Herald1,herald2,herald-boden} the third  condition holds if the
set
of perturbations is {\em abundant}.

Next,  let $\cP$ denote the set of $C^r$ functions (for some fixed large
$r$),
$f:SU(n)^N\to \rr$, invariant under the conjugation action of
$SU(n)$, i.e. $f(rg_1r^{-1},\cdots,rg_Nr^{-1})=f(g_1,\cdots,g_N)$.

Now let $P|_{D^2}$ denote the restriction of the principal bundle $P$ to
the disc $D^2=\gamma_i(D^2\times\{1\})$. Given a smooth connection $A\in
\cA$, consider the map
$\Phi_A:P|_{D^2}\to SU(n)^N$ defined  by the rule
$\Phi_A(p)=(g_1,\cdots ,g_N)$ where $pg_i\in P_x$ is the endpoint of
the
parallel lift of  the loop $\gamma_i(x,u), \ u\in [0,2\pi]$  with
respect
to
$A$ starting at
$p$ (where $x\in D^2$ is the image of $p$ via the bundle projection).
Then
$\Phi_A(pr)=(r^{-1}g_1 r,\cdots,r^{-1}g_N r)$ for $r\in SU(n)$. It
follows
that if $f\in \cP$, $f\circ \Phi_A$ descends to a well-defined function
on
$D^2$.

If the principal bundle is trivialized, $P\cong X\times SU(n)$, then one
can view parallel lifting as the holonomy, a function from loops to
$SU(n)$, and in our context this allows us to define
$$\Hol:\cA\times D^2\to SU(n)^N, \ \ (A,x)\mapsto
(\Hol_{\gamma_1(x,-)}(A),\cdots \Hol_{\gamma_N(x,-)}(A)).$$
     Thus we denote the function $f\circ \Phi_A$  by
$x\mapsto f(\Hol(A,x))$. Notice  that $f(\Hol(A,x))$ is independent of
the
trivialization, and is unchanged if we change $A$ by a gauge
transformation.

Fix a smooth positive cut--off function $\eta$ on $D^2$ which vanishes
near
the boundary. The space of {\em admissible perturbations} is defined to
be
the space of functions of the form
\begin{equation}\label{eq2.5}h_f:\cA\to \rr, \ A\mapsto
\int_{D^2}f(\Hol(A,x))\eta(x)\ d^2x.\end{equation}
     Thus we view  the admissible perturbation functions as being
parameterized by the (Fr{\'e}chet) space $\cP$.

Taubes and Floer show that $h_f$ extends to a smooth function on
$\cA_2$.  Note that $d(cs +h_f)_A= dcs_A + d({h_f})_A $ and
$\hess(cs+h_f)_A=\hess(cs)_A+\hess(h_f)_A$.  Thus one can define
$$Q_{A,f}:T_A\cA\cong\Om^1(X;\ad P)\to T_A\cA $$ by  the rule
$$\langle
Q_{A,f}(b),c\rangle=\hess(h_f)_A(b,c)=\tfrac{\partial^2}{\partial t 
\partial s}|_{s=t=0}\ h_f(A+sb +tc).$$  Then setting $d_{A,f}:=
d_A-4\pi^2
*Q_{A,f}$ we see that $-\tfrac{1}{4\pi^2} *d_{A,f}$
     is the operator defined by the Hessian of $cs+h_f$ at $A$.

A connection $A$ is called {\em $f$-perturbed flat} if it is a critical
point of $cs +h_f$  or, equivalently,  if $*F(A)=4\pi^2
\grad(h_f)_A$. As in Section \ref{sec2.3}, we form the sequence
\begin{equation}\label{eq2.6}0\to \Om^0(X;\ad P)\xrightarrow{d_A}
\Om^1(X;\ad P)\xrightarrow{d_{A,f}}
     \Om^2(X;\ad P)\xrightarrow{d_A}
      \Om^3(X;\ad P)\to 0.\end{equation} The sequence \eqref{eq2.6} is a
complex if $A$ is
$f$-perturbed flat, and for any connection $A$ and perturbation $f$ the
operator obtained from  \eqref{eq2.6} (compare to \eqref{foldedupop}),
\begin{equation}
\begin{split} &D_{A,f}:\Om^0(X;\ad P)\oplus \Om^1(X;\ad P) \to
\Om^0(X;\ad
P)\oplus
\Om^1(X;\ad P),\\
&D_{A,f}(b,c) =(d_A^* c ,*d_{A,f} c  + d_A b )
\end{split}
\end{equation}
is a self--adjoint operator called the {\em perturbed
twisted odd signature operator}.

The operator $*Q_{A,f}$ is independent of the choice of
Riemannian metric on $X$ since the definition of
$\hess(h_f)_A(b,c)$ does not use the metric, and since
$$\hess(h_f)_A(b,c)=\langle Q_{A,f}(b),c\rangle=-\int_X\tr(c\wedge
*Q_{A,f}(b)).
$$ It follows that $d_{A,f}$ is metric independent. In particular, if
$A$
is $f$-perturbed flat, then the cohomology of \eqref{eq2.6}, and hence
the
dimension of the kernel of $D_{A,f}$, is independent of the metric on
$X$.

We will give a proof  in Lemma \ref{estimate}  of Taubes' observation
that for any
$A$ and $f$, the operator
$D_{A,f}$ is Fredholm and in fact a compact perturbation of
$D_A$. Assuming this,  we conclude the following.

\begin{lemma}\label{independence} Suppose that $f_0$ and $f_1$ are
perturbations,
$A_0$ and
$A_1$ are connections so that $A_0$ is $f_0$ perturbed flat and
$A_1$ is
$f_1$ perturbed flat.  If $f_t$ is a path of perturbations from $f_0$ to
$f_1$ and
$A_t$ is a path of connections from $A_0$ to $A_1$ then
the spectral flow $\SF(D_{A_t,f_t})$ is independent of the choice of
paths
$A_t$ and
$f_t$ and independent of the  choice of Riemannian metric.
\end{lemma}
\begin{proof}
Since the space of connections and the space of perturbations are
contractible and the spectral flow is a homotopy invariant,
$\SF(D_{A_t,f_t})$ is independent of the choice of paths $A_t$ from
$A_0$ to $A_1$ and path $f_t$ from $f_0$ to $f_1$. Since the dimensions
of
the kernels of $D_{A_0,f_0}$ and $D_{A_1,f_1}$ are independent of the
Riemannian metric, it similarly follows that   $\SF(D_{A_t,f_t})$ is
independent of the choice of Riemannian metric.
\end{proof}

\subsection{Locality of $Q_{A,f}$}\label{secloc}
     The operator $D_{A,f}$ is not a differential operator, but it has
the
property that it differs from a differential operator by an operator
which
is  ``localized'' in the neighborhood of the solid tori used to
construct the perturbations.

Write
$$D_{A,f}(b,c)= D_A(b,c) + (0, -4\pi^2 Q_{A,f}(c))$$ so that $D_{A,f}$
differs from the Dirac operator $D_A$ by
$-4\pi^2 Q_{A,f}$.

\begin{prop}\label{lem2.1} Let $S\subset X$ denote the union of the
images
of the embeddings $\gamma_i$, $S=\cup_i\gamma_i(D^2\times S^1)$.
\begin{enumerate}
\item If $b,b'\in \Om^1(X;\ad P)$ and $b|_{S}=b'|_{S}$, then
$Q_{A,f}(b)=Q_{A,f}(b')$.
\item For any $b\in\Om^1(X;\ad P)$, $Q_{A,f}(b)$ vanishes outside of
$S$.
\item If $A,A'$ are connections and $A|_S=A'|_S$, then
$Q_{A,f}(b)=Q_{A',f}(b)$ for all $b\in\Om^1(X;\ad P)$.
\end{enumerate}

\end{prop}
\begin{proof} Suppose that  $b\in\Om^1(X;\ad P)$ vanishes on $S$.  Then
for each
$c\in\Om^1(X;\ad P)$,
$$\langle Q_{A,f}(b),c\rangle=\hess(h_f)_A(b,c)=
\tfrac{\partial^2}{\partial s  \partial
t}|_{s=t=0}\int_{D^2}f(\Hol(A+sb+tc,x))\eta(x)d^2x=0,$$
since the parallel lifts of the loops $\gamma_i(x,u), \ u\in[0,2\pi]$
with
respect to the connection $A+sb+tc$ are independent of $s$.  Hence
$Q_{A,f}(b)=0$. Part (1) now follows since $Q_{A,f}$ is linear and
$b-b'$
vanishes on $S$.

For part (2), Let $b\in\Om^1(X;\ad P)$ be arbitrary. Let $U\subset X\setminus S$
     and choose $c\in\Om^1(X;\ad P)$ supported in $U$.  Then
     $$\langle Q_{A,f}(b),c\rangle=\langle b,Q_{A,f}(c)\rangle=0$$
     by the same argument given for part (1) since $c$ vanishes on
     $S$.  Since the $L^2$ inner product of $Q_{A,f}(b)$ with $c$
     vanishes for every $c$ supported on $U$, $Q_{A,f}(b)$ must vanish on
     $U$.

Part (3) is clear.

\end{proof}

     Proposition \ref{lem2.1} implies that  the
restriction of
$D_{A,f}$  to
     a codimension 0 submanifold whose boundary misses
$S=\cup_i$Image$(\gamma_i)$ is well--defined. To be precise, suppose
that  $Y\subset X$ is a compact
     codimension 0 submanifold. We consider two cases:
$S\subset$interior$(Y)$
     and $S\subset X\setminus Y$.

     If $S\subset$interior$(Y)$, $A$ is any connection on $Y$ and
       $c$ is a 1-form  on $Y$ with values in $\ad
     P|_Y$, extend $A$ and $c$ arbitrarily to $\tilde{A}, \tilde{c}$ on
$X$.
Then
     $Q_{\tilde{A},f}(\tilde{c})$ vanishes outside of $S$ and is
independent
of the choice of extensions. Thus $Q_{A,f}$ is well-defined as an
operator
on $Y$.
     Since $D_A$ is a differential operator it is locally defined and
hence
has a well defined restriction to $Y$. It follows that
$D_{A,f}=D_A-4\pi^2
Q_{A,f}$ is well defined on $Y$.

If $Y\subset X\setminus S$, then a similar argument shows that $D_{A,f}$
restricts
to the differential operator $D_A$ on $Y$.

\subsection{An estimate}  We next give an  argument to show that
$D_{A,f}$ differs from $D_A$ by an operator which is bounded on $L^2$.
This result is essentially contained in Taubes' article \cite{taubes}.
See
also
\cite[Proposition 2.8]{herald-boden}.

\begin{lemma} \label{estimate} There exists a constant
$C$   so that for all $A\in
\cA, f\in\cP$, and $b\in \Om^1(X;\ad P)$
$$\|Q_{A,f}(b)\|_{L^2}\leq C \|f\|_{C^2}\|b\|_{L^2}.$$ In particular,
$D_{A,f}$ is a compact perturbation of
$D_A$, viewed as operators $L^2_1\to L^2$.
\end{lemma}
\begin{proof}

Denote $\Hol(A+sb+tc,x)$ by $H_x(s,t)$,  so
$H_x:\rr^2\to  SU(n)^N$, and let
$f:SU(n)^N\to \rr$ be invariant under the conjugation action of $SU(n)$.
Then
$$
\hess(h_f)_A (b,c) = \left.\tfrac{\partial^2}{\partial
      s\partial t}\right|_{s=t=0} \int_{D^2} f(H_x(s,t))\eta(x) d^2x =
\int_{D^2}
\left.\tfrac{\partial^2}{\partial
      s\partial t}\right|_{s=t=0}  f(H_x(s,t)) \eta(x)d^2x.
$$

We may think of $H_x:\rr^2\to SU(n)^N\subset
\cc^{n^2N}$ and  extend  $f$  to a smooth function
$\cc^{n^2N}\to \rr$.    By the chain and the product rule  for vector
valued functions  we have
\begin{eqnarray*}
\lefteqn{\left.\tfrac{\partial^2}{\partial
      s\partial t}\right|_{s=t=0}  f(H_x(s,t))}\\ & = & D_{H_x(0,0)} f
\circ
\left.\tfrac{\partial^2}{\partial
      s\partial t}\right|_{s=t=0} H_x(s,t) +
      \left(\left.\tfrac{\partial}{\partial s}\right|_{s=0}
      H_x(s,0)\right)^T \circ D_{H_x(0,0)}^2 f \circ
\left.\tfrac{\partial}{\partial t}\right|_{t=0}
      H_x(0,t)
\end{eqnarray*}

Then
\begin{eqnarray*}
\lefteqn{\left|\left.\tfrac{\partial^2}{\partial
      s\partial t}\right|_{s=t=0}  f(H_x(s,t))\right|}\\
& \leq & \|D_{H_x(0,0)} f \| \cdot \|
\left.\tfrac{\partial^2}{\partial
      s\partial t}\right|_{s=t=0} H_x(s,t) \|+
      \|\left.\tfrac{\partial}{\partial s}\right|_{s=0}
      H_x(s,0)\| \cdot \|D_{H_x(0,0)}^2 f\| \cdot \|
\left.\tfrac{\partial}{\partial t}\right|_{t=0}
      H_x(0,t)\|
\end{eqnarray*}

Given  $c\in \Om^1(X;\ad P)$, $i\in \{1,\cdots ,N\}$,  and
$x\in D^2$,  let
$\tilde c_{i,x}:[0,1]\to su(n)$ be defined as follows.  Trivialize $\ad
P$
over the disc $D^2\subset X$, so $\ad P|_{D^2}\cong D^2\times su(n)$.
Parallel  translation  by $A$ along the path $g_{i,x}:[0,1]\to X,
t\mapsto\gamma_i(x,e^{2\pi i t})$ trivializes $g_{i,x}^*(\ad P)$ over
$[0,1]$, and this way the pullback $g_{i,x}^*(c)$ is viewed as a 1-form
with
$su(n)$ coefficients. Contracting with the unit speed vector field
$\tfrac{\del}{\del t}$ yields the function
$\tilde c_{i,x}:[0,1]\to su(n)$, i.e.~
$\tilde c_{i,x}=g_{i,x}^*(c)(\tfrac{\del}{\del t})$.

We have the estimate
$$\int_{D^2} \int_0^1 |\tilde c_{i,x}(\nu)|\eta(x)\ d\nu\ d^2x
\leq \int_{D^2} \int_0^1 |c(\nu)|\eta(x) \ d\nu\ d^2x\leq C_1\|
c\|_{L^1(S^1\times D^2)}\leq C\| c\|_{L^1(X)}$$ for some constants $C_1,
C$ independent of $c$.

Given a  connection $B$ let
$h_{i,x}(B)$ denote the holonomy of $B$ around the loop
$\gamma_i(x,-)$. Thus
$$H_x(s,t)=(h_{1,x}(A+s b+t c),
\cdots,h_{1,x}(A+s b+t c)).$$   Lemma 2.6    of
\cite{herald-boden} establishes the formulas:
\begin{eqnarray*}
\left.\tfrac{d}{dt}\right|_{t=0} h_{i,x}(A+t b) & = & h_{i,x}(A)
\int_0^1
\tilde b_{i,x}(\nu) d\nu\\
\left.\tfrac{\partial^2}{\partial s\ \partial t}\right|_{s=t=0}
h_{i,x}(A+s b+tc) & = &
h_{i,x}(A)
\int_0^1\int_0^\nu (\tilde b_{i,x}(\nu)\tilde c_{i,x}(\mu) + \tilde
c_{i,x}(\nu)\tilde b_{i,x}(\mu)) d\mu d\nu.
\end{eqnarray*}

Notice that $h_{i,x}(A)\in SU(n)$, and so is bounded independently of
$A$.  Thus  (letting $C$ denote possibly different constants)
$$
\int_{D^2}  |  \left.\tfrac{d}{dt}\right|_{t=0} H_x(s,0) | \ \eta(x)\
d^2x
\leq\sum_{i=1}^N C\|b\|_{L^1(X)} = NC\|b\|_{L^1(X)}\leq
C\|b\|_{L^2(X)}.$$
The last inequality follows from Cauchy-Schwarz.

Similarly
\begin{eqnarray*}
     \lefteqn{\int_{D^2} | \left.\tfrac{\partial^2}{\partial
      s\partial t}\right|_{s=t=0} H_x(s,t)  | \ \eta(x) \ d^2x }\\
& \leq  & C \int_{D^2} \int_0^1 |b(\nu)|d\nu \int_0^1|c(\mu)|d\mu d^2x
   \leq C \left\|\int_0^1 |b(\nu)|d\nu\right\|_{L^2(D^2)} \cdot
\left\|\int_0^1|c(\mu)|d\mu\right\|_{L^2(D^2)}\\
& \leq &  C \| b\|_{L^2(S^1 \times D^2)} \cdot \| c \|_{L^2(S^1 \times
D^2)}
\leq C
\|b\|_{L^2(X)}
\cdot \|c\|_{L^2(X) }.
\end{eqnarray*}

Also
$$ \|D_{H_x(0,0)} f\| + \|D^2_{H_x(0,0)} f\| \leq C \|f\|_{C^2}.
$$
Thus for any $b,c\in \Om^1(X;\ad P)$,
$$  \bigl|\hess(h_f)_A (b,c)\bigr|  \leq C \|f\|_{C^2}
\|b\|_{L^2} \|c\|_{L^2},
$$
and hence $\|Q_{A,f}\|\le C \|f\|_{C^2}$.

For the  last statement of Lemma \ref{estimate}, the difference
$D_{A,f}-D_A:L^2_1\to L^2$ equals the restriction  of
$-4\pi^2Q_{A,f}$ to $L^2_1$. This restriction  can  be expressed as the
composite of the compact inclusion
$L^2_1\subset L^2$ and the bounded operator $-4\pi^2Q_{A,f}$, and hence
is
compact since the compact operators form a 2-sided ideal.
\end{proof}

\section{The Poisson operator and the \Cald\ for perturbed Dirac
operators
on manifolds with boundary} A compact manifold with non-empty boundary
always embeds  as a codimension 0 submanifold of a closed manifold
(e.g.~
its double) in such a way that a given principal bundle and $SU(n)$
connection extend. It follows from Proposition \ref{lem2.1} that given:
\begin{enumerate}
\item a 3-manifold with boundary $X$,
\item a principal bundle $P\to X$,
\item a connection $A$ on $P$,
\item a collection of embeddings $\gamma_i:D^2\times S^1\to
$interior$(X)$, and
\item a function $f\in \cP$
\end{enumerate} one can unambiguously define operators $Q_{A,f}$ and
$D_{A,f}$ on
$X$.  In this context we will assume that the connection $A$ is in
cylindrical form near the boundary. This means that in a collar
neighborhood
$[0,\epsilon)\times \partial X$  of $\partial X$, $A$ is the pullback
of a
connection $a$ on $\partial X$ with respect to the projection
$[0,\epsilon)\times \partial X\to \partial X$ to the second factor. Then
the operator $D_{A,f}$ equals $D_A$ on the collar, and this implies that
$D_{A,f}$ takes the Atiyah--Patodi--Singer form (for details see
\cite[Definitions 2.2 and 2.3]{BHKK})
\begin{equation} D_{A,f}|_{[0,\epsilon)\times \partial
X}=D_{A}|_{[0,\epsilon)\times
\partial X}=
\gamma(\tfrac{\partial}{\partial u} + S_a). \label{G4.1}
\end{equation}

Here $S_a:\Om^*(\partial X;\ad P)\to \Om^*(\partial X;\ad P)$ is the
tangential operator, given by the formula
\begin{equation}\label{eq4.45}S_a(\alpha_0,\alpha_1,\alpha_2)=(*d_a\al_1
,
-*d_a\al_0-d_a*\al_2, d_a*\al_1)
\end{equation}
    and $\gamma$ is the bundle isomorphism given by the formula
\begin{equation}\label{eq4.4}
\gamma(\alpha_0,\alpha_1,\alpha_2)
=(-*\al_2,*\al_1,*\al_0)\end{equation}
(in the  formulas \eqref{eq4.4} and \eqref{eq4.45}, $*$ refers to the
Hodge
$*$ operator for the metric on $\partial X$). The bundle isomorphism
$\gamma$ satisfies $\gamma^2=-\Id$.

\subsection{The invertible double}

So far we have shown that $D_{A,f}$ is the sum of a Dirac operator $D_A$
and  an  $L^2$--bounded operator $-4\pi^2 Q_{A,f}$. The operator
$Q_{A,f}$
is not a differential operator but it has certain locality properties as
explained in Section \ref{secloc}. On a manifold with boundary we
arranged
so that $D_{A,f}$ has a product form in a collar of the  boundary.

In the study of Dirac operators on manifolds with boundary, their
\Cald\ and their Poisson operator play a central role (cf.
\cite{booss-w}).
However, since $Q_{A,f}$ is far from being a differential operator,
existence of the \Cald\ and of the Poisson operator is a priori not
clear.
We therefore treat the existence question and related issues. For the
convenience of the reader our presentation will be more or less
self--contained.

Let $X$ be a compact connected Riemannian manifold with boundary
$\partial
X=\Sigma$ and $D$ a self--adjoint Dirac type operator on the Hermitian
vector bundle
$E\to X$. We assume that we are in a product situation, i.e. that there
is
a collar $U=[0,\epsilon)\times \Sigma$ of $\Sigma$ in which $D$ takes
the
form
\begin{equation}
      D=\gamma\bigl(\frac{d}{dx}+B\bigr),
\label{G4.10}
\end{equation} where $\gamma$ is a bundle endomorphism and $B$ is a
first
order self--adjoint elliptic differential operator on $E|\Sigma$. One
has
the fundamental relations
\begin{equation}
     \gamma^2=-\Id,\quad \gamma^*=-\gamma,\quad\text{and}\quad
     B\gamma=-\gamma B.\label{G4.11}
\end{equation} Note that $D_A$ is of this form \eqref{G4.1}.

To model $Q_{A,f}$ we assume abstractly that we are given an auxiliary
self--adjoint bounded linear operator $V$ in $L^2(X,E)$. In most cases
$V$
will be a pseudo--differential operator of order $0$, but this is not
assumed here. More importantly, we assume the following locality
property
(cf. Sec. \ref{secloc}).
\begin{equation}
\text{There is a compact domain
$S\subset X\setminus [0,\epsilon)\times \Sigma$ such that $M_\varphi
V=0$
for all $\varphi\in C^\infty_0(X\setminus S)$.}
\label{pseudolocality}
\end{equation}
Here $M_\varphi$ denotes the operator of multiplication by $\varphi$.
Note that by self--adjointness we automatically also have
\begin{equation}
      VM_\varphi=0\quad\text{for}\quad \varphi\in C^\infty_0(X\setminus
S).
      \label{G4.13}
\end{equation} Moreover, if $f,g\in L^2(E)$ such that $f|W=g|W$ in a
neighborhood $W$ of $S$ then choosing $\psi\in C^\infty_0(W)$ with
$\psi\equiv 1$ in a neighborhood of $S$ we find
\begin{equation} V(f-g)=V((1-\psi)(f-g))=(VM_{1-\psi})(f-g)=0.
\label{G4.14}
\end{equation}

The property \eqref{pseudolocality} implies that the restriction of
$D+V$
to the collar $U$ equals $D$ and hence takes the  product form
\eqref{G4.10}. This and self--adjointness is all that is required to
prove
Green's formula (see \cite{booss-w})
\begin{equation}
\langle (D+V) f,g\rangle_{L^2(X,E)}-
\langle f, (D+V) g\rangle_{L^2(X,E)}=\langle f|_\Sigma,
\gamma(g|_\Sigma)\rangle_{L^2(\Sigma,E)}.
\end{equation}

An important tool for studying elliptic boundary problems is the {\em
invertible double construction} and the associated Poisson operator.
Since
it will be important in the sequel we recall the main facts (cf.
\cite{booss-w}).

Let $\tilde X:=X_+\cup_{\Sigma} X_-$ where $X_\pm=X$. Let
$\tilde E\to \tilde X$ be the vector bundle obtained by gluing
$E$ with gluing function $\gamma$ over $\Sigma$. More precisely,
consider
the bundle $E\times\{\pm\} \to X\times\{\pm\}$. Then
$$\tilde{X} =X\times\{\pm\}/\sim,  \text{ where } (x,+)\sim (y,-) \text{
if and only if }
    x=y\text{ and }
    x\in\Sigma.$$  We write $X_\pm:=X\times\{\pm\}$.  Also,
$$\tilde{E}=E\times\{\pm\}/\sim\text{ where } (x,v,+)\sim (y,w,-)
\text{ if
and only if }
    x=y\in\Sigma \text{ and }  \gamma(v)=w.
     $$  Denote by $m:\tilde X\to\tilde X$ the reflection interchanging
$X_+$
and $X_-$, i.e.  $m(x,\pm)=(x,\mp)$. Then $m$ is covered by a
bundle
map
$\tilde{m}:\tilde{E}\to \tilde{E}$ defined by
$\tilde{m}:(x,v,\pm)=(x,\mp
v,\mp)$.  Notice that the restriction of $\tilde{m}$ to $
\tilde{E}|\Sigma$
equals $\gamma$ since if $x\in \Sigma,$
$\tilde{m}:(x,v,\pm)=(x,\mp v, \mp)\sim (x, \gamma(v),\pm).$ The map
$\tilde{m}$
preserves the Hermitian metric on $E$  and
$\tilde{m}^2=-\Id$.

The maps $m$ and  $\tilde{m}$  induce a unitary map
\begin{equation}
\Phi:L^2(\tilde X,\tilde E)\to L^2(\tilde X,\tilde E),\quad
\Phi f(x):=\tilde{m} ( f(m(x)))
\end{equation}
satisfying
$\Phi^2=-\Id$.  Thus $\Phi$ is defined by the diagram
\[\begin{diagram}
\node{\tilde{E}}\node{\tilde{E}}\arrow{w,t}{\tilde{m}}\\
\node{\tilde{X}}\arrow{n,r}{\Phi
f}\arrow{e,t}{m}\node{\tilde{X}}\arrow{n,r}{f}
\end{diagram}\]

Then $D$ extends naturally to a Dirac operator $\tilde D$
on
$\tilde E$ satisfying
\begin{equation}
       \Phi^* \tilde D\Phi=-\tilde D.\label{G4.15}
\end{equation} Similarly, in view of \eqref{pseudolocality}, $V$ extends
to a bounded operator $\tilde V$ satisfying $\Phi^*\tilde V\Phi=-\tilde
V$
which has the same locality properties as $V$ with
$\tilde S=S\cup m(S)$ instead of $S$.

Alternatively, we may view $\tilde D$ as a self--adjoint realization of
a
well--posed boundary value problem of $D\oplus (-D)$ acting on the
bundle
$E\oplus E$ over $X$. Namely, putting
\begin{equation}
          \begin{split}
              \Psi:&L^2(\tilde X,\tilde E)\longrightarrow L^2(X,E\oplus
E)\\
                   &\Psi f:= f|X_+\oplus (\Phi f)|X_+
          \end{split}\label{G4.16}
\end{equation} we find
\begin{equation}
        \Psi\tilde D\Psi^{-1}=D\oplus (-D)
\end{equation} and
\begin{equation}
      \Psi(\operatorname{dom}\tilde D)=\Psi(L^2_1(\tilde E))=
        \bigl\{f\in L^2_1(X,E\oplus E)\,\big|\, (f|\Sigma)_2=\gamma
(f|\Sigma)_1\bigr\}.\label{G4.17}
\end{equation} It is indeed not difficult to show that
$(f|\Sigma)_2=\gamma (f|\Sigma)_1$ is a well--posed boundary condition
for
$D\oplus (-D)$. The boundary operator of this boundary condition,
an orthogonal projection in the pseudodifferential
Grassmannian (see \cite{booss-w}),
is given by
\begin{equation}
       \frac 12\begin{pmatrix} \Id & \gamma\\ -\gamma &
\Id\end{pmatrix}.\label{G4.18}
\end{equation}

We introduce the usual notations for restriction, extension, and trace
operators:

\begin{defn}  Let $r^\pm$ be restriction of sections on $\tilde X$ to
$X_\pm$ and
$e^\pm$ extension by $0$ from $X_\pm$. $\varrho^\pm,\varrho$ denote the
trace maps
\begin{equation}\begin{split}
       &L^2_s(X_\pm,E)\longrightarrow L^2_{s-1/2}(\Sigma,E),\\
       &L^2_s(\tilde X,\tilde E)\longrightarrow L^2_{s-1/2}(\Sigma,E),
                   \end{split}\quad s>1/2.\label{G4.19}
\end{equation} Similarly, for $0<|t|<\epsilon$ we denote by $\varrho_t$
the trace map which restricts sections to $\{t\}\times\Sigma$. For
$\varrho_t$
\eqref{G4.19} holds accordingly.

Note that, by the locality property \eqref{pseudolocality} of $V$, we
have
$r^+\tilde V=Vr^+.$
\end{defn}

For convenience we abbreviate $T:=D+V$ resp. $\tilde T:=\tilde D+\tilde
V$. The pseudolocality property of $\tilde V$ implies  that $\tilde T$
has
well--defined restrictions to $X_{\pm}$.  The following fact is
well--known:

\begin{lemma}\label{ML-S4.3} For arbitrary $s\in [0,\infty)$ the trace
map
$\varrho^+$ extends to a bounded linear map
\[   \ker T\cap L^2_s(X,E)\to L^2_{s-1/2}(\Sigma,E).\] More precisely,
if
$f\in \ker T\cap L^2_s(X,E)$ then
$\varrho^+(f)=\lim\limits_{t\to 0+}\varrho_t(f)$ exists in
$L^2_{s-1/2}(\Sigma,E)$ and the so defined $\varrho^+$ is bounded.
\end{lemma}
With some care the lemma could even be stated for all $s\in\mathbb{R}$.
However,
Sobolev spaces of negative order on manifolds with boundary are a
nuisance.
Since we will not need them we content ourselves to the case $s\ge 0$.
\begin{proof} For $V=0$ the Lemma is well--known
(\cite[Thm. 13.1]{booss-w}, cf. also \cite[Chap. I]{lions-magenes}).
If $Tf=0$ then $Df=0$ in a
collar of $\Sigma$. Since the result is local in a collar of
$\Sigma$, we reach the conclusion.\end{proof}

In light of Lemma \ref{ML-S4.3} we can safely introduce  the following
notation.

\begin{defn} Set
\begin{align*}
          Z_\pm^s&:=\bigl\{ f\in L^2_s(X_\pm,E)\,\big|\, \tilde
T|_{X_\pm}
f=0\bigr\},\quad s\ge 0,\\
          \La_\pm^s&:=\varrho^\pm Z_\pm^{s+1/2},\quad s\ge -1/2,\\
          Z_0&:=\bigl\{f\in L^2_1(\tilde X,\tilde E)\,\big|\, \tilde T
f=0\bigr\}.
\end{align*}
\end{defn} Note that from elliptic regularity we immediately conclude
that
if
$f\in L^2(\tilde X,\tilde E)$ and $\tilde Tf=0$ then
$f\in L^2_1(\tilde X,\tilde E)$. Since $V$ is only assumed to be
$L^2$--bounded, this cannot be improved. Also by elliptic theory,
$\tilde T$ is a Fredholm operator (it is a relatively compact
perturbation
of $\tilde D$) and hence $\dim Z_0<\infty$.

\begin{prop}\label{ML-S4.5} If $f\in Z_0$ then $f$ vanishes in a collar
of
$\Sigma$ and hence $Z_0=r^+ Z_0\oplus r^-Z_0$.

In particular, $\tilde D$ is an invertible operator,   called the
\emph{invertible double of $D$}.
\end{prop}
\begin{proof} We apply Green's formula and find
\begin{equation}\begin{split}
        \| \varrho f\|^2&=-\langle\varrho^+r^+f,\gamma \gamma
\varrho^+r^+f\rangle= -\langle\varrho^+r^+f,\gamma \varrho^+r^+(\Phi
f)\rangle\\
             &= -\langle Tf,\Phi f\rangle_{L^2(X_+,E)}+\langle f,T\Phi
f\rangle_{L^2(X_+,E)}=0
\end{split}
\end{equation} since $f\in Z_0$.

Thus we have $\varrho f=0$. Since the Dirac operator has the weak unique
continuation property (cf. \cite{booss-ucp}) this implies that $f$
vanishes
in a collar neighborhood of $\Sigma$. If $V=0$ then the weak unique
continuation property and the connectedness of $X$ imply together with
the
proven vanishing statement for $f\in Z_0$ that $\tilde D$ is
invertible.\end{proof}

\begin{defn} We say that $D+V$ has the {\em weak unique continuation
property with respect to $\partial X$} if $f=0$ is the only element of
$Z_+^{1/2}$ with $\varrho^+f=0$.
\end{defn}

In general we cannot expect $D+V$ to have the weak unique continuation
property. However, for small $V$ it is true:

\begin{prop}\label{WUCP} The operator $D+V$ has the weak unique
continuation property with respect to $\partial X$ if and only if
$\tilde
D+\tilde V$ is invertible. This is in particular the case if  $\|\tilde
D^{-1}\tilde V\|<1$.
\end{prop}
\begin{remark} Booss--Bavnbek, Marcolli, and Wang \cite{BMW} prove
the weak unique continuation property for a class of perturbed Dirac
operators
(perturbations may be nonlinear)
arising in Seiberg--Witten theory. However, our perturbations
are in general not admissible in their sense (cf. \cite[Def. 2.5]{BMW}).
An operator $\mathfrak{P}$ is admissible in the sense of \cite{BMW} if
there is a locally bounded function $P$ such that
$|\mathfrak{P}u(x)| \le P(u,x)|u(x)|$. This implies in particular that
$\mathfrak{P}u(x)=0$ if $u(x)=0$, i.e.
$\mathfrak{P}$ is a local operator. If $\mathfrak{P}$ is a linear
operator
on smooth sections then this locality automatically implies that
$\mathfrak{P}$
is a differential operator of order $0$ (i.e. a bundle endomorphism). So
the point in \cite{BMW} is that $\mathfrak{P}$ is allowed to be
nonlinear.
\end{remark}

\noindent{\it Proof of Proposition \ref{WUCP}.} In view of \eqref{G4.15}
and Prop.
\ref{ML-S4.5} the weak unique continuation
property with respect to $\partial X$ implies $Z_0=0$ and hence the
invertibility of $\tilde D+\tilde V$.

Conversely, assume that $\tilde D+\tilde V$ is invertible and consider
$f\in Z_+^{1/2}$ with $\varrho^+f=0$. Then $e^+f\in Z_0=\{0\}$ and hence
$f=0$.

    $\|\tilde D^{-1}\tilde V\|<1$ implies the invertibility of
$\tilde D+\tilde V$ by means of the Neumann series.
\qed

\subsection{Existence of Poisson operator and \Cald}\label{sec41}

Let $P_{Z_0}$ be the orthogonal projection onto $Z_0$ and denote by
$\varrho^*$ the $L^2$--dual of $\varrho$. In other words, for a
distributional section $f$ of $\tilde E|\Sigma$ and a test function
$\varphi\in C^\infty(\tilde X,\tilde E)$ one has $(\varrho^*
f,\varphi):=(f,\varrho
\varphi)$.
In view of Proposition \ref{ML-S4.5} one
infers from this formula immediately that
$P_{Z_0}\varrho^*=0$.

Furthermore, denote by $\tilde G$ be the pseudoinverse of $\tilde T$,
i.e.
\begin{equation}
      \tilde Gf:=(\Id-P_{Z_0})(\tilde T+P_{Z_0})^{-1}f=
      \begin{cases} \tilde T^{-1}f,& f\in Z_0^\perp,\\
                            0,& f\in Z_0.
      \end{cases}\label{G4.20}
\end{equation} Note that since $\tilde T$ maps $L^2_s\to L^2_{s-1}, 0\le
s\le 1$ the pseudoinverse maps $L^2_s\to L^2_{s+1}, -1\le s\le 0$.
Together with the mapping properties of the trace map we infer that
$\tilde G\varrho^*$ maps $L^2_s(\Sigma,E)$ continuously into
$L^2_{s+1/2}(\tilde X,\tilde E)$ for each $-1/2\le s<0$. It is a
fundamental fact that for $r^\pm \tilde D^{-1}\varrho^*$ this can be
improved. Namely, one has:

\begin{thm}\label{Poissonop-a} $r^\pm\tilde D^{-1}\varrho^*$ maps
$L^2_s(\Sigma,E)$ continuously into $L^2_{s+1/2}(X_\pm,E)$ for all
$s\ge -1/2$.
\end{thm}

For a proof see \cite{booss-w}. It will also follow from our discussion
of
the continuous dependence of the \Cald\ below. First we want to show
that
the previous theorem easily extends to hold for $\tilde T$ and hence a
\Cald\ can be constructed for $\tilde T$.

\begin{thm}\label{Poissonop-b}
$r^\pm\tilde G\varrho^*$ maps $L^2_s(\Sigma,E)$ continuously into
$Z_\pm^{s+1/2}$ for $-1/2\le s\le 1/2$.

Furthermore, we have the following resolvent identity relating $\tilde
G$
and $\tilde D^{-1}$:
\[\tilde G=\tilde D^{-1}-\tilde G \tilde V\tilde D^{-1}-P_{Z_0}\tilde
D^{-1}.
\]
\end{thm}
\begin{proof} The resolvent identity follows from
\begin{equation}\begin{split}
       \tilde G\tilde V\tilde D^{-1}= \tilde G(\tilde T-\tilde D)\tilde
D^{-1}\\
                 = (\Id-P_{Z_0})\tilde D^{-1}-\tilde G.
\end{split}\end{equation} The operator $P_{Z_0}\tilde D^{-1}$ is
bounded
from $L^2_s(\tilde X,\tilde E)$ to  $L^2_1(\tilde X,\tilde E)$ for $s\ge
-1$  (see    Prop. \ref{ML-S4.5} and the comment preceding it).
Moreover,
$\tilde G\tilde V\tilde D^{-1}$ maps
$L^2_{-1}(\tilde X,\tilde E)$ continuously into $L^2_1(\tilde X,\tilde
E)$.
   From Theorem \ref{Poissonop-a} and the resolvent identity
we thus conclude that
$r^\pm\tilde G\varrho^*$ maps $L^2_s(\Sigma,E)$ continuously into
$L^2_{s+1/2}(X_\pm,E)$.

Moreover, since
$\tilde T \tilde G\varrho^*=(\Id-P_{Z_0})\varrho^*=\varrho^*$ and
$\varrho^*$ is supported on $\Sigma$, we find that the operator actually
maps into $Z_\pm^{s+1/2}$.
\end{proof}

The restriction  $-1/2\le s\le 1/2$ in Theorem \ref{Poissonop-b}
is due to the low regularity
assumptions on $V$. If $V$ is a pseudodifferential operator, then
Theorem \ref{Poissonop-a} holds verbatim for $\tilde T$ instead of
$\tilde D$.
The previous proof can easily be modified to this case.

With Theorem \ref{Poissonop-b} and Lemma  \ref{ML-S4.3}  in place we can
make the following definitions.

\begin{defn}\begin{enumerate}
\item Define $$K_\pm:=K(T)_\pm:=\pm r^\pm\tilde
G\varrho^*\gamma:L^2_s(\Sigma,E)\to Z^{s+1/2}_{\pm}, \ -1/2\leq s\leq
1/2.$$
Then $K:=K_+$ is called the {\em Poisson operator} for
$T=D+V$.
\item Define the  projection $$C_\pm:=C_\pm(T):=\varrho^\pm
K(T)_\pm:L^2_s(\Sigma,E)\to
L^2_s(\Sigma,E), \ -1/2\leq s\leq 1/2.$$   Then $C:=C_+$ is called the {\em
\Cald\ }of
$T$.
\item Taking $s=0$, the subspace $$\La :=\La(T):=
\La^0_+(T) =\operatorname{range}C_+(T)\subset
L^2(\Sigma,E)
$$ is called the {\em Cauchy data space} of $T$.
\end{enumerate}
\end{defn}

\begin{prop}\label{ML-S4.9} $C_\pm(D+V)$ are complementary orthogonal
projections in $L^2(\Sigma,E)$ with
$\operatorname{range}C_\pm=\La_\pm^0$.

Furthermore, $\gamma C_\pm\gamma^*=C_\mp$. In other words the subspaces
$\La_\pm^0$ are Lagrangian subspaces with respect to the symplectic form
$\omega(x,y)=\langle x,\gamma y\rangle_{L^2}.$
\end{prop}

Recall that a Lagrangian subspace is a subspace $V\subset L^2(\Sigma,E)$
satisfying $\gamma V=V^\perp$.

\begin{proof}
$K_\pm$ maps $L^2(\Sigma,E)$ into $Z_\pm^{1/2}$ by Theorem
\ref{Poissonop-b}. Hence in view of Lemma \ref{ML-S4.3} we find that
$C_\pm$ is
$L^2$--bounded with
$\operatorname{range} C_\pm\subset \La_\pm^0$.

We show  $\La_+^0\perp \La_-^0$: pick $\xi\in \La_+^0,
\eta\in \La_-^0$ and choose $f\in Z_+^{1/2}, g\in Z_-^{1/2}$ with
$\xi=\varrho^+f, \eta=\varrho^-g$. Then Green's formula yields
\begin{equation}
       \begin{split}
(\xi,\eta)&=(\varrho^+f,\varrho^-g)=-(\varrho^+f,\gamma\varrho^+(\Phi
g))\\
         &=-(T f,\Phi g)_{L^2(X_+,E)}+(f,T\Phi g)_{L^2(X_+,E)}=0.
       \end{split}
\end{equation} Hence $\La_+^0\perp \La_-^0$.

Next we show that $C_++C_-=\Id$: pick $\xi\in L^2(\Sigma,E)$ and $f\in
C^\infty(\tilde X,\tilde E)$.   Green's formula gives
\begin{equation}\begin{split}
            ((C_++C_-)\xi,\gamma\varrho f)_\Sigma
&= (\varrho^+r^+\tilde
G\varrho^*\gamma \xi,\gamma\varrho f)_\Sigma-(\varrho^-r^-\tilde
G\varrho^*\gamma\xi,\gamma\varrho f)_\Sigma\\
      &=(TK_+\xi,f)_{X_+}-(\tilde G\varrho^*\gamma\xi,Tf)_{X_+}+
        (TK_-\xi,f)_{X_-}-(\tilde G\varrho^*\gamma\xi,\tilde Tf)_{X_-}\\
      &=-(\tilde G\varrho^*\gamma\xi,\tilde
Tf)_X=-((\Id-P_{Z_0})\varrho^*\gamma
\xi,f)_X\\
      &=(-\gamma\xi,\varrho f)_\Sigma=(\xi,\gamma\varrho f)_\Sigma.
           \end{split}
\end{equation}
Since $\gamma\varrho: C^\infty(\tilde X,\tilde E)\to C^\infty (\Sigma;
E)$ is surjective, $(C_++C_-)\xi=\xi$.

In sum, $C_\pm$ are $L^2$--bounded with orthogonal range $\subset
\La_\pm^0$
and $C_++C_-=\Id$. Hence $C_\pm$ must be the orthogonal projections onto
$\La_\pm^0$.

Finally, consider $\xi=\varrho^+f\in \La_+^0, f\in Z_+^{1/2}$. Then, by
\eqref{G4.15} we have $\Phi f\in Z_-^{1/2}$ and thus
$\gamma \xi=\varrho^-(\Phi f)\in \La_-^0$. Hence $\gamma$ interchanges
$\La_+^0$ and $\La_-^0$. Since $\La_+^0\perp \La_-^0$ the spaces
$\La_\pm^0$ must
be
Lagrangian and we are done.
\end{proof}

\subsection{Continuous dependence of the \Cald}\label{sec42} Let us
first
describe the set--up in which we want to study the dependence of
the \Cald\ on its input data:  We consider the class $\cD$ of Dirac
operators on $E$ which are in cylindrical form near $\Sigma$ and the
class
$\mathcal{V}$ of auxiliary operators satisfying
\textnormal{\eqref{pseudolocality}}. We say that a family $D_t\in\cD$
varies continuously if all coefficients of $D_t$ (in any local chart)
vary continuously. Moreover, we equip
$\cV$ with the usual norm topology of $\cB(L^2(X,E))$.
Finally, we put $(\cD\times \cV)(n)
:=\bigl\{(D,V)\in\cD\times\mathcal{V}\,\big|\, \dim Z_0( D+V)=n\bigr\}$.

We start with some elementary statements on the continuous dependence of
$P_{Z_0}$ and $\tilde G$:

\begin{prop}\label{contGtilde} The maps
\begin{equation}
\begin{split}P_{Z_0}, \tilde G&:(\cD\times\cV)(n)\to \cB(L^2(\tilde
X,\tilde E),
L^2_1(\tilde X,\tilde E)),\\
     &P_{Z_0}(D,V):=\proj(\ker(\tilde D+\tilde V)),\\
     &\tilde G(D,V):=(\Id-P_{Z_0}(D,V))(\tilde D+\tilde
V+P_{Z_0}(D,V))^{-1},
   \end{split}\label{ML-G3.21}
\end{equation}
are continuous.
\end{prop}
Note that this statement is slightly stronger than the continuity
of $P_{Z_0}, \tilde G$ as maps into $\cB(L^2(\tilde X,\tilde E))$.

\begin{proof} Operators in $\cD$ are in cylindrical form near $\Sigma$,
hence
$\cD\times \cV\to \cB(L^2_1(\tilde X,\tilde E),L^2(\tilde X,\tilde E)),
(D,V)\mapsto \tilde
D+\tilde V$ is continuous. Thus also
$\mathbb{C}\times\cD\times \cV\to \cB(L^2_1(\tilde X,\tilde
E),L^2(\tilde
X,\tilde E)), (\lambda,D,V)\mapsto \tilde
D+\tilde V-\lambda$ is continuous. For a fixed pair $(D_0,V_0)\in
(\cD\times\cV)$ there is an $\varepsilon>0$ such that $0$ is the only
spectral point of $\tilde D_0+\tilde V_0$ in the closed disc
$\{z\in\mathbb{C}\,|\,
|z|\le \varepsilon\}$. By continuity this remains true for $(D,V)$ in
a sufficiently small neighborhood of $(D_0,V_0)$ in
$(\cD\times\cV)(n)$. For these $(D,V)$ we have
\begin{equation}
      P_{Z_0}(D,V)=\frac{1}{2\pi i}\oint_{|z|=\varepsilon}(\lambda-\tilde
      D-\tilde V)^{-1}d\lambda.
\end{equation}
Since the inversion is continuous this proves the continuity statement
on
$P_{Z_0}$.

The continuity statement on $\tilde G(D,V)$ now follows immediately from
the defining formula \eqref{ML-G3.21}.
\end{proof}

The resolvent identity in Theorem
\ref{Poissonop-b} shows that for $(D,V)\in(\cD\times \cV)(n)$
the Poisson operator and hence $C_+(D+V)$
depend continuously on $V$: namely $K_+(D+V)=K_+(D)-r^+(\tilde
G(D,V)\tilde
V\tilde D^{-1}+P_{Z_0}(D,V)\tilde D^{-1})\varrho^*\gamma$ and by the
previous Proposition
$$(D,V)\mapsto r^+(\tilde G(D,V)\tilde V\tilde D^{-1}
+P_{Z_0}(D,V)\tilde D^{-1})\varrho^*\gamma$$
is continuous
$(\cD\times \cV)(n)\to \cB(L^2_{-1/2}(\Sigma,E),L^2_1(X_+,E))$.
However, the dependence on
$D$ itself is more subtle. In order to explore
this we have to look a bit closer at the construction of the \Cald.

Fix a $c\in \mathbb{R}$ with $c\not\in\operatorname{spec} B$.
We abbreviate the spectral projections
$P_{>c}:=1_{(c,\infty)}(B),P_{\le c}:=1_{(-\infty,c]}(B)$ etc. Then we
set,
for $\xi\in L^2_s(\Sigma,E)$,
\begin{equation}
    Q(B)\xi(x):=e^{-xB}P_{\ge c}\xi, \quad x\ge 0.
\end{equation}

We first need to study the mapping properties of $Q(B)$ with respect
to Sobolev norms. The results we are going to present are fairly
standard
applications of the spectral theorem. However, to be as self--contained
as possible we   present the details.

Let  $\epsilon >0$, $\epsilon <1$ and fix a smooth cut--off
function $\varphi:\rr\to [0,1]$ with
$\varphi(x)=1$ for $ |x| \le\epsilon/2$  and $\varphi(x)=0$ for $|x|\ge
\epsilon$.
Let $\psi:\rr\to \rr$ be any function in $
C^\infty_0((-\epsilon,\epsilon)\setminus\{0\})$,  for example
$\psi=\varphi'$.
In the following proposition we only  use the restrictions of $\varphi$
and $\psi$ to
$\{x\ge 0\}$.

\begin{prop}\label{ML-S3.11}\hfill

\textnormal{1. } For all
$s,s'\in\mathbb{R}$ we have $e^+\psi
Q(B)\in\cB(L^2_s(\Sigma,E),L^2_{s'}(\tilde X,\tilde E))$.

\textnormal{2.}
For $-1/2\le s\le 1/2$ the operator
$\varphi Q(B)$ maps
$L^2_s(\Sigma,E)$ continuously into $L^2_{s+1/2}(X_+,\tilde E)$.
\end{prop}
\begin{proof} First we need to recall the definition of Sobolev spaces
on a  manifold with boundary.  Since $B$ is elliptic
the Sobolev norms on $L^2_s(\Sigma,E)$ can be defined using
$\Id+|B|$:
\begin{equation}
       \langle f,g\rangle_{L^2_s}:=\langle (\Id+|B|)^s f,(\Id+|B|)^s
g\rangle_{L^2}
\end{equation}
for $f,g\in L^2_s(\Sigma,E)$. Furthermore, for a nonnegative integer
$k\in \mathbb{Z}_+$ a Sobolev norm for $L^2_k(\mathbb{R}_+\times
\Sigma,E)$
is given by
\begin{equation}
      \|f\|_{L^2_k}^2=\sum_{j=0}^k \int_0^\infty
\|\partial_x^j(\Id+|B|)^{k-j}f(x)\|_{L^2(\Sigma,E)}^2dx.\label{Sobnorm}
\end{equation}
The Sobolev norms on $X_+$ are obtained by  patching together  the
usual interior Sobolev norms with the norms \eqref{Sobnorm}.

1. By complex interpolation it suffices to prove the claim for $s'\in
\mathbb{Z}_+$. Applying the Leibniz rule to
$\partial^j(\Id+|B|)^{k-j}\psi
Q(B)$ we see that in view of \eqref{Sobnorm} it is in fact sufficient
to show that for any $k\in\mathbb{Z}_+$ the operator $(\Id+|B|)^k\psi
Q(B)$ is bounded $L^2_s(\Sigma,E)\to L^2(X_+,E)$.

Furthermore, for $\xi\in L^2_s(\Sigma,E)$
we have $(\Id+|B|)^{s}\xi\in L^2(\Sigma,E)$
and 
\[(\Id+|B|)^k\psi Q(B)\xi=(\Id+|B|)^{k-s}\psi
Q(B)((\Id+|B|)^{s}\xi)\]
so that  it suffices to prove the $L^2$--boundedness of
$(\Id+|B|)^k\psi Q(B)$ for all $k$.

Next, since $\psi\in C^\infty_0(0,\epsilon)$ there is an $x_0>0$
such that $\psi(x)=0$ for $x\le x_0$. For $x_0\le x\le\epsilon$ the
spectral
theorem gives
\begin{equation}
     \|(\Id+|B|)^ke^{-xB}P_{\ge c}\|_{L^2(\Sigma,E)}
=\sup_{\lambda\ge c}|(1+\lambda)^ke^{-x\lambda}|\le C(k,x_0,c,\epsilon)
\end{equation}
and thus for $\xi\in L^2(\Sigma,E)$
\begin{equation}\begin{split}
       \|\psi (\Id+|B|)^kQ(B)\xi\|_{L^2(X_+,E)}^2&\le \|\psi\|_\infty^2
\int_{x_0}^\epsilon
              \|(\Id+|B|)^ke^{-xB}P_{\ge c}\xi\|_{L^2(\Sigma,E)}^2dx\\
       &\le \epsilon \|\psi\|^2 C(k,x_0,c,\epsilon)^2 \|\xi\|^2.
                \end{split}
\end{equation}
  The proof of 1. now follows since $\psi$ has compact support in $(0,\infty)$
and hence in this case extension by $0$ 
is bounded $L^2_{s'}(X_+)\to L^2_{s'}(\tilde{X})$ for any $s'$.

\smallskip

2. Let $(E_\lambda)_{\lambda\in\mathbb{R}}$ be the spectral resolution
of
$B$. Then we estimate for $\xi\in L^2(\Sigma,E)$
\begin{equation}\label{ML-basicestimate}
    \begin{split}
        \|\varphi&(\Id+|B|)^{1/2}Q(B)\xi\|
\le \int_0^\infty \| (\Id+|B|)^{1/2}\varphi(x)e^{-xB} P_{\ge
c}(B)\xi\|^2dx\\
     &\le \int_c^\infty \int_0^\epsilon (1+|\lambda|) e^{-2x\lambda}dx
d\langle
        E_\lambda \xi,\xi\rangle\\
      &\le C(c,\epsilon) \|\xi\|^2.
    \end{split}
\end{equation}

For $\xi\in L^2_{-1/2}(\Sigma,E)$ we can now estimate using
\eqref{ML-basicestimate}
\begin{equation}\begin{split}
      \|\varphi Q(B)\xi\|_{L^2(X_+,E)}&=\|\varphi
      (\Id+|B|)^{1/2}Q(B)(\Id+|B|)^{-1/2}\xi\|_{L^2(X_+,E)}\\
        &\le \|\varphi (\Id+|B|)^{1/2}Q(B)\|_{L^2(\Sigma,E)\to
      L^2(X_+,E)}\|(\Id+|B|)^{-1/2}\xi\|_{L^2(\Sigma,E)}
                \end{split}
\end{equation}
proving the continuity of $\varphi Q(B)$ from $L^2_{-1/2}(\Sigma,E)$
to $L^2(X_+,E)$.

Next consider $\xi\in L^2_{1/2}(\Sigma,E)$. Then analogously we find
(taking
$\psi=\varphi'$ and applying part 1)
\begin{equation}\begin{split}
      \|\partial_x\varphi Q(B)\xi\|_{L^2(X_+,E)}&\le \|\varphi'
      Q(B)\xi\|_{L^2(X_+,E)}+ \|\varphi
      B(\Id+|B|)^{-1/2}Q(B)(\Id+|B|)^{1/2}\xi\|_{L^2}\\
     &\le C_1 \|\xi\|_{L^2}+C_2\|(\Id+|B|)^{1/2}\xi\|_{L^2}
                \end{split}
\end{equation}
and
\begin{equation}\begin{split}
      \|\varphi (\Id+|B|)Q(B)\xi\|_{L^2(X_+,E)}
         &= \|\varphi (\Id+|B|)^{1/2}Q(B)(\Id+|B|)^{1/2}\xi\|_{L^2}\\
     &\le C(c,\epsilon)\|(\Id+|B|)^{1/2}\xi\|_{L^2}.
                \end{split}
\end{equation}
In view of \eqref{Sobnorm} this proves the continuity of $\varphi Q(B)$
from $L^2_{1/2}(\Sigma,E)$ to $L^2_1(X_+,E)$.

The complex interpolation method \cite[Sec. 4.2]{taylor} now yields the
assertion for
all $s\in [-1/2,1/2]$.
\end{proof}

\newcommand{\Ell}{\operatorname{Ell}}
\begin{prop}\label{ML-continuity} Let $\Ell^1(\Sigma,E)$ be the set
of first order self--adjoint elliptic differential operators on
$E|\Sigma$.
Then for $\psi\in C^\infty_0((-\epsilon,\epsilon)\setminus{0})$ the map
\[
    \Ell^1(\Sigma,E)_c:=\{B\in\Ell^1(\Sigma,E)\,|\,
c\not\in\operatorname{spec} B\}
   \longrightarrow \cB(L^2(\Sigma,E),L^2(X_+,E)),\quad B\mapsto \psi Q(B)
\]
is continuous.
\end{prop}
\begin{proof} Again let $x_0>0$ be such that $\psi(x)=0$ for $0\leq x
\le x_0$.
By the spectral theorem the map
\begin{equation}
     [x_0,\epsilon]\times \Ell^1(\Sigma,E)_c\longrightarrow
\cB(L^2(\Sigma,E)),
    \quad (x,B)\mapsto e^{-xB} P_{\ge c}(B)
\end{equation}
is continuous.

Now fix $B_0\in \Ell^1(\Sigma,E)_c$ and a $\delta>0$. By compactness
there is a neighborhood $W$ of $B_0$ such that for $
B\in W$ and all $x_0\le x\le \epsilon$
\begin{equation}
     \|e^{-xB}P_{\ge c}(B)-e ^{-xB_0}P_{\ge c}(B_0)\|_{L^2(\Sigma,E)}\le
     \delta/\sqrt{\epsilon}.
\end{equation}
Thus for $\xi\in L^2(\Sigma,E)$
\begin{equation}
    \|\psi Q(B)\xi-\psi Q(B_0)\xi\|^2\le C^2\int_{x_0}^\epsilon
     \|e^{-xB}P_{\ge c}(B)-e ^{-xB_0}P_{\ge c}(B_0)\|_{L^2(\Sigma,E)}^2dx
      \|\xi\|^2
    \le C^2\delta^2\|\xi\|^2.
\end{equation} 
Thus $  \|\psi Q(B) -\psi Q(B_0) \| \le C\delta$,    completing the
proof.
\end{proof}

Let
\begin{equation}\label{defapproxCald}
       R(B)\xi(x):=\begin{cases}
                  -e^{-xB}P_{\ge c}\gamma \xi  & \text{for }x\ge 0, \\
                   e^{-xB}P_{\le c}\gamma \xi  & \text{for } x<0,
     \end{cases}
\end{equation}
and define the {\em approximate    Poisson operator
}   $K^0_\pm:=\pm r^\pm \varphi R(B)\gamma$.

Applying Proposition \ref{ML-S3.11} and \ref{ML-continuity}
(and its obvious counterparts for $x\le 0$) to $R(B)$ gives:
\begin{prop}\label{approxCald}\textnormal{1.} If $\psi\in
C^\infty_0((-\epsilon,\epsilon)\setminus\{0\})$ then for all
$s,s'\in\mathbb{R}$,
$$\psi
R(B)\in\cB(L^2_s(\Sigma,E),L^2_{s'}(\tilde X,\tilde E)).$$  Moreover,
as an element of $\cB(L^2(\Sigma,E),L^2(\tilde X,\tilde E))$,
$\psi R(B)$ depends continuously on the coefficients of $B$ as long
as $c\not\in\operatorname{spec} B$.

\textnormal{2.} The map  $\varphi R(B)$ maps
$L^2(\Sigma,E)$ continuously into
$L^2_{1/2}(X_-\amalg X_+,\tilde E)$. 
\end{prop}

If $f\in C^\infty(\tilde X,\tilde E)$ then Green's formula gives
\begin{equation}\begin{split}
       (\varphi R\xi,\tilde T f)_{\tilde X}&=(r^+\varphi R\xi,r^+\tilde
Tf)_{X^+}+(r^-\varphi R\xi,r^-\tilde
       Tf)_{X^-}\\
       &= (-\gamma P_{\ge c}\gamma\xi,\varrho f)_\Sigma-(\gamma P_{\le
       c}\gamma\xi,\varrho f)_\Sigma+(S(B)\xi,f)_{\tilde X}\\
       &=((\varrho^*+S(B))\xi,f)_{\tilde X}
       \end{split}
\end{equation}
with $S(B)\xi=\varphi'\gamma R\xi$. In view of Prop.
\ref{approxCald} we have in particular that $S(B)\in
\cB(L^2(\Sigma,E),L^2(\tilde X,\tilde E))$ and that it depends
continuously on the coefficients of $B$. Note furthermore, that
$S(B)$ maps in fact into
$L^2_{s,\textrm{comp}}(\tilde X\setminus\Sigma,\tilde E)$ for all $s\ge
0$. Hence we
have
shown that
\begin{equation}
      \tilde T\varphi R=\varrho^*+S(B).\label{ML-G4.21}
\end{equation}

The definition of $K_\pm^0$ immediately gives
\begin{equation}
         \begin{split}
           \varrho^+K_+^0\xi&:=\lim_{t\to 0+}\varrho_tK_+^0\xi=P_{\ge
c}\xi,\\
           \varrho^-K_-^0\xi&:=\lim_{t\to 0-}\varrho_tK_-^0\xi=P_{\le
c}\xi.
         \end{split}\label{ML-G4.22}
\end{equation} Note that we assumed $c\not\in\operatorname{spec} B$.

\eqref{ML-G4.21} and \eqref{ML-G4.22} allow to express
$K_\pm(\tilde D+\tilde V)$ and $C_\pm(\tilde D+\tilde V)$ in terms of
the
operator $R(B)$. Namely, we have in view of \eqref{ML-G4.21}
\begin{equation}
        (\Id-P_{Z_0})\varphi R=\tilde G\tilde T\varphi R=
           \tilde G\varrho^*+\tilde G S(B)
\end{equation} and hence
\begin{equation}\begin{split}
       K_+ &=r^+\tilde G\varrho^*\gamma\\
           &=r^+(\Id-P_{Z_0})\varphi R\gamma-r^+\tilde GS(B)\gamma\\
           &=K_+^0-r^+P_{Z_0}\varphi R\gamma-r^+\tilde GS(B)\gamma.
          \end{split}
\end{equation} Similarly,
\begin{equation}
       K_-=K_-^0+r^-P_{Z_0}\varphi R\gamma+r^-\tilde GS(B)\gamma.
\end{equation} Finally, this gives     the formulas:
\begin{equation}\begin{split}
       C_+&=P_{\ge c}(B)-\varrho^+\tilde G S(B)\gamma,\\
       C_-&=P_{\le c}(B)+\varrho^-\tilde G S(B)\gamma.
\end{split}\label{ML-G4.26}
\end{equation}
In view of Proposition \ref{contGtilde} and Proposition \ref{approxCald}
the operator
$\varrho^\pm\tilde G S(B)\gamma$ maps $L^2(\Sigma,E)$ continuously into
$L^2_{1/2}(\Sigma,E)$ and, as long as $\dim Z_0$ is constant, it
depends continuously on $D$ and $V$.

Since $c\not\in\operatorname{spec}(B)$ the operators $P_{\ge c}(B),
P_{\le
c}(B)$ depend continuously on $B$. Note that the choice of $c$ is
irrelevant. To show the continuity of
$C_+$ at $D$ and $V$ we choose a $c$ not in $\operatorname{spec}(B)$.
Then
we will have $c\not\in\operatorname{spec}(B)$ in a full neighborhood of
$D$ and
$V$. Of course, $S(B)$ also depends on $c$.

Summing up we have proved:

\begin{thm}\label{continuityCald} We consider the class $\cD$ of Dirac
operators on $E$ which are in cylindrical form near $\Sigma$ and the
class
$\mathcal{V}$ of auxiliary operators satisfying
\textnormal{\eqref{pseudolocality}}. We say that a family $D_t\in\cD$
varies continuously if all coefficients of $D_t$ (in any local chart)
vary
continuously. Moreover, we equip
$\mathcal{V}$ with the usual norm topology of $\cB(L^2(X,E))$. Then for
each
$n\in\mathbb{Z}_+$ the map
\[ \bigl\{(D,V)\in\cD\times\mathcal{V}\,\big|\, \dim Z_0( D+V)=n\bigr\}
      \longrightarrow \cB(L^2(\Sigma,E)), \quad (D,V)\mapsto C( D+ V)
\] is continuous.
\end{thm}

Note also that in view of \eqref{ML-G4.26} the operator $C-P_{\ge
0}(B)$
maps $L^2$ continuously into $L^2_{1/2}$, i.e. `is  an operator of order
$-1/2$'. If $V$ is a pseudodifferential operator then one easily shows
that $\varrho^\pm \tilde G S(B)$ is a smoothing operator and we recover,
though with a different argument, the well--known fact that $C-P_{\ge
0}(B)$ is smoothing
\cite[Prop. 2.2]{Scott}, \cite[Prop. 4.1]{Grubb}.

\subsection{Application to the perturbed Dirac operator}

We turn to our specific operator $D_{A,f}=D_A-4\pi^2 Q_{A,f}$.   The
discussion in Section \ref{secloc} shows that
$Q_{A,f}$ satisfies \eqref{pseudolocality}. Hence the general results of
Sections \ref{sec41} and \ref{sec42} apply. In particular there exist
the
Poisson operator $K_+(D_{A,f})$ and the Calder{\'o}n projector
$P_{A,f}:=C_+(D_{A,f})$. Next we show that in a neighborhood of the flat
connections and for ``small'' $f$  the weak unique continuation property
with respect to $\partial X$ holds for $D_{A,f}$.

\begin{lemma}\label{nbdU}  There is an open neighborhood
$U\subset
\cA$ of the set of flat connections and a constant $C$ (which depends
only
on $X$ as a Riemannian manifold) so that given any
$A\in U$,
$$\|\tilde D_{A}^{-1}\|\leq C.$$
     \end{lemma}
     \begin{proof} Given a flat connection
$A$ and a gauge transformation $g$ over $X$ the operators $D_A$ and
$D_{g\cdot A}$ are unitarily equivalent (i.e. $D_{gA}(a)=gD_A
g^{-1}(a)$),
and this equivalence extends to $\tilde D_A$ since $\cG$ acts on the Lie
algebra coefficients and $\gamma=\pm *$ acts on the differential forms,
thus $g_x\gamma_x(v)=\gamma_x(g_x(v))$ for $x\in \Sigma$. Hence the
spectra of
$\tilde{D}_A$ and $\tilde{D}_{gA}$ coincide.

The space $\cM$ of flat connections modulo gauge  transformations on $X$
is compact, being homeomorphic to the  compact  real algebraic variety
$\Hom(\pi_1X,SU(n))/\text{conjugation}$. Thus   there is a lower bound,
say $k$, on the absolute value of the smallest eigenvalue of
$\tilde{D}_A$ as $A$ varies in the set of flat connections on $X$.  It
follows that the largest eigenvalue of $\tilde{D}_A^{-1}$ is bounded by
$1/k$, and hence $\|\tilde{D}_A^{-1}\|\leq 1/k$. Set
$C=2/k$, then by continuity the set $U=\{A\in\cA\ |\
     \|\tilde{D}_A^{-1}\|\leq C\}$ is a neighborhood of the set of flat
connections.
\end{proof}

\begin{prop}\label{lem4.2} There is an $\epsilon>0$ so that given any
connection $A$ in the neighborhood $U$ of Lemma \textnormal{\ref{nbdU}},
and any
     $f\in\cP$ satisfying $\|f\|_{C^r}<\epsilon$, then
$\tilde D_A-4\pi^2\tilde Q_{A,f}$ is invertible and hence has the weak
unique continuation property with respect to
$\partial X$.
\end{prop}
\begin{proof}  The Neumann series shows that $\tilde D_A-4\pi^2\tilde
Q_{A,f}$ is invertible  if $\|\tilde D_A^{-1} 4\pi^2 \tilde
Q_{A,f}\|<1$.

By Lemmas \ref{nbdU} and
\ref{estimate} there exists a constant $C$ independent of $A$ and $f$
such that
$$\|\tilde D_A^{-1}4\pi^2 \tilde Q_{A,f}\|\le 4\pi^2 \|\tilde D_A^{-1}\|
\|\tilde Q_{A,f}\|\le C\|f\|_{C^2}.$$ Hence we may choose
$\epsilon:=\tfrac{1}{C}$. Proposition \ref{WUCP} then sows that
$D_{A,f}$
has the weak unique continuation property.
\end{proof}

\begin{thm}\label{thm5.3} The map
$U\times B_{\cP}(0,\epsilon)\to \cB(L^2(\Sigma,E))$ sending a pair
$(A,f)$ to the Calder\'on projector $C_{A,f}$ of $D_{A,f}$ is
continuous,
and hence the Cauchy data spaces $\Lambda_{A,f}:=
\operatorname{range}(C_{A,f})$ vary continuously (in the gap metric on
the
Grassmannian of Lagrangian subspaces of $L^2(\Sigma,E)$) with respect to
$A,f$.
\end{thm}
\begin{proof} By Prop. \ref{lem4.2} we have $Z_0(\tilde D_A-4\pi^2\tilde
Q_{A,f})=0$ for all $(A,f)\in \cA\times B_{\cP}(0,\epsilon)$. Thus the
claim follows from Theorem \ref{continuityCald}.
\end{proof}

\section{An application}\label{application} Theorem \ref{thm5.3} allows
us
to use the technology developed in \cite{nico} and \cite{daniel-kirk}
to
study how the spectral flow of a path of operators
$D_{A_t,f_t}$ on a closed 3-manifold behaves with respect to a
decomposition along a separating surface. In particular, the  proof of
Nicolaesu's adiabatic limit theorem and Lemma 3.2 of \cite{daniel-kirk},
which consider the effect of stretching the collar of $X$ on the Cauchy
data spaces,  do not depend on the operators being Dirac operators, but
only on the fact that they have the Atiyah--Patodi--Singer form on a
collar.  Thus     Theorem 5.1 of
\cite{daniel-kirk}, and all its consequences (Section 6 of
\cite{daniel-kirk}) are true for the operators
$D_{A,f}$.

This will be used crucially in forthcoming   calculations of the
$SU(3)$ Casson invariant for Brieskorn spheres
\cite{BHK-pqr}.  We finish this article with a useful  application of
our
main result. We outline a different proof of a result of Taubes which is
the main ingredient in his proof that his invariant equals Casson's
$SU(2)$
invariant for homology 3-spheres, and is the starting point for proofs
of
various approaches to problems known as the  Atiyah--Floer conjecture.
It is also  one of the simplest
non-trivial examples to which this machinery can be applied, since
this situation avoids the delicate problem of singularities in the
moduli space.

\bigskip

   First, we set up  some notation. Let $X$ be
a homology 3-sphere.  Let
$P\to X$ be  the trivial principal
$SU(2)$ bundle  and let $\cA(X)$ denote the ($L^2_1$ completion of the)
space of
$SU(2)$ connections on
$P$. Denote by $\cG(X)$ the ($L^2_2$ completion of the) gauge group
$\Map(X,SU(2))$.  Then let $\cA^*(X)$ be the subspace of $\cA(X)$
consisting of irreducible connections, that is, those connections whose
stabilizer in  $\cG(X)$ is equal to  the center $\{\pm \Id \}
\subset\cG(X)$. Thus $\cA^*(X)/\cG(X)$ is a smooth Hilbert manifold.

    Suppose that
$X$ is given a Heegard decomposition
$X=X_1\cup_\Sigma X_2$ of genus $g>2$. Choose    embeddings $\gamma_i$
of
solid tori so that $S=\cup_i\gamma_i(D^2\times S^1)$   is contained in
the
interior of
$X_2$, and  choose   $f\in \cP$ so   that  the restriction of $cs+h_f$
to
the space $\cA^*(X)$ of irreducible connections  has only
non-degenerate
critical points (i.e.~the kernel of
$D_{A,f}$ is zero at every irreducible perturbed flat connection $A$).
In
this context the fact  that such $f$ are dense in
   any  neighborhood of $0\in\cP$   is proven in
\cite{taubes} and
\cite{floer}. We will  restrict $f$ further below.

It follows that the critical points of
$cs+h_f:\cA^*(X)/\cG(X)\to \rr/\zz$ are finite and isolated.
    Let $\cM_f^*(X)$ denote the critical set, i.e.~ the moduli space of
irreducible perturbed flat connections on $X$:
$$ \cM_f^*(X)=\{A\in\cA^*(X)\ | \ d(cs+h_f)(A)=0\}/\cG(X).$$

   To a pair $[A_0], [A_1]\in
\cM^*_f(X)$, Taubes assigns the sign
$$\Sign_T([A_0],[A_1])= (-1)^{\SF(D_{A_t,f})} $$ where $A_t$ is any path
in $\cA^*(X)$ joining lifts of $[A_0]$ and $[A_1]$. In fact an
application
of the  index theorem shows that the mod $8$ reduction of
$\SF(D_{A_t,f})$
depends only on the gauge equivalence classes $[A_0]$ and $[A_1]$.

Restricting irreducible connections to the pieces of the  Heegard
decomposition leads to a diagram (with the obvious notation)
\[
\begin{diagram}\dgsquash[3/4]
\node[2]{\cM^*(X_1)}\arrow{se}\\
\node{\cM^*_f(X)}\arrow{ne}\arrow{se}\node[2]{\cM^*(\Sigma)}\\
\node[2]{\cM^*_f(X_2)}\arrow{ne}
\end{diagram}
\] All arrows are embeddings. The space $\cM^*(\Sigma)$ is a smooth
     $6g-6$ dimensional  symplectic manifold. (We will describe the
symplectic structure below.)    The subspaces
$\cM^*(X_1)$ and
$\cM_f^*(X_2)$ are smooth  and
orientable $3g-3$ dimensional Lagrangian  submanifolds
(\cite{taubes,Herald1}). They intersect transversely in the compact
$0$--dimensional manifold $\cM_f^*(X)$ for   appropriate $f$,
namely those perturbations $f$ so that $cs+h_f$ has only
non-degenerate
critical points in $\cA^*(X)$.   Notice that since
$S$ lies in the interior of
$X_2$,  perturbed flat connections restrict to flat connections on $X_1$
and $\Sigma$.

We also assume that $\|f\|_{C^2}$ is small enough so  the following
requirements  hold:
\begin{enumerate}
\item The space $\cM^*_f(X_2)$ is diffeomorphic to the flat (i.e.
unperturbed) moduli space $
\cM^*(X_2)$. This is possible since, as $f$ varies, $\cM^*_f(X_2)$ varies
by
a (Legendrian) cobordism (\cite{Herald1}) and so for small enough $f$
the
cobordism is a product. (Some care must be taken since $\cM^*_f(X_2)$
is non-compact,
but the results of \cite{Herald1}  can be used to reach this conclusion.
Alternatively, Taubes shows that the perturbations can be chosen to fix
the reducible
stratum.)
\item The perturbation should be small enough so that $\cM^*_f(X_2)$
lies
in the neighborhood $U$ specified by Lemma \ref{nbdU}. That this  is
possible follows again from  \cite{Herald1}.
\item The perturbation should be small enough so that  the requirement
$\|f\|<\epsilon$ of  Proposition \ref{lem4.2} holds on $X_2$.
\end{enumerate}

Then Theorem
\ref{thm5.3} implies that the map taking $A\in \cM^*_f(X_2)$ to the
Cauchy data space $\Lambda_{A,f}$ for $D_{A,f}$ is continuous.

    In this language Casson assigns the sign (after  orienting
    the  moduli spaces $\cM^*(X_1),\ \cM_f^*(X_2)$ and $\cM^*(\Sigma)$)
\[\Sign_C([A_0],[A_1])=\begin{cases} 1 & \text{if } 
 \cM^*(X_1)\cdot_{[A_1]}\cM_f^*(X_2)=\cM^*(X_1)\cdot_{[A_0]}\cM_f^*(X_2)
 \text{ and,}\\
 -1  &  \text{otherwise}, 
\end{cases}\]
where
$\cM^*(X_1)\cdot_{[A_1]}\cM_f^*(X_2)\in \{\pm 1\} $ denotes the local
transverse
intersection  number of $\cM^*(X_1)$ and $\cM_f^*(X_2)$ in
$\cM^*(\Sigma)$
at
$[A_1]$.   

Taubes' proof that his invariant equals Casson's invariant (up to an
overall sign which depends on a method for fixing particular
orientations
of all the moduli spaces)  reduces to the assertion that
$\Sign_T=\Sign_C$. Notice that we are looking at relative signs, so that
the  validity of $\Sign_T=\Sign_C$ is independent of the choice of
orientations of the moduli  spaces.

We begin with a lemma whose proof appears  to be well  known to those
who
know it well  (e.g. \cite{dostoglou-salomon}), although  we were  unable
to  find  an explicit statement in the literature.  The
     statement is implicit in \cite{daska},  but is described in the
language of stable holomorphic bundles.  We are indebted to G.
Daskalopoulos for help with  the argument.

\begin{lemma}\label{lemma6.1} The  space $\cF^*(\Sigma)$ of irreducible
flat
$SU(2)$ connections over a closed 2-manifold $\Sigma$ of genus greater
than 2 is   simply connected. The moduli space
$\cM^*(\Sigma)$ of such connections is also  simply connected.

\end{lemma}

\begin{proof}   Before we start the proof, we remark that in the
statement
of  Lemma \ref{lemma6.1}, the space $\cF^*(\Sigma)$ can mean either the
smooth flat irreducible connections, or the $L^2_r$ completion of this
space for any $r$ large enough so that the
$L^2_{r+1}$ completion of the gauge group consists of continuous gauge
transformations,   e.g.~$r>0$.

We first set up notation.  We consider the trivial rank 2 smooth
vector bundle
$E=\cc^2\times\Sigma\to \Sigma$  endowed with a Hermitian
metric.  We fix a holomorphic structure on
$\Sigma$.  The space $\cA$ of
$SU(2)$ connections on $E$  is  
identified with the space  of holomorphic structures with
holomorphically trivial determinant on
$E$ via the standard construction (see
\cite{Atiyah-Bott}) since $\Sigma$ is 1 complex--dimensional and so every
structure is integrable. We denote by   $\cA^s$  the subspace of stable
holomorphic structures on $E$ (see \cite{Atiyah-Bott}).  We take the
$L^2_r$ completions of these spaces.

Denote by $\cF \subset \cA$ the subspace of flat connections, and
by $\cF^* $ the subspace of irreducible flat connections.   We
  take $\cG$ to be the ($L^2_{r+1}$
completion of the) $SU(2)$ gauge group
$\text{Map}(\Sigma,SU(2))$.

The slice theorem (see e.g.~\cite{donaldson-kronheimer}) shows that
the
quotient map $\cA\to
\cA /\cG $, when restricted to the irreducible connections
$\cA^*$    gives a locally trivial map, and hence a fibration
$\cA^*\to \cA^*/\cG$ with fiber
$\cG/\pm\Id$. Restricting to the flat irreducible connections
gives the  fibration
\begin{equation}\label{eq6.8}\cG /{\pm\Id}\to
\cF^*\to \cM^*.
\end{equation}

That $\cM^*$ is simply connected follows immediately from  \cite[Theorem
7.1]{daska} and the fact
 that  $\cM^*$ is
homeomorphic to
the space of stable rank 2 holomorphic bundles over $\Sigma$ with
trivial
determinant (see \cite{narasimhan-sheshadri, Atiyah-Bott}).

The calculations in \cite{daska} show that
$\pi_2(\cM^*)=\zz\oplus\zz/2$ and it is not  hard to show  that
$\pi_1(\cG/{\pm\Id})=\zz\oplus\zz/2$.  If we knew that   the
map induced by the connecting homomorphism  in the long exact sequence
of
homotopy groups for the fibration \eqref{eq6.8} were an isomorphism (or
merely surjective), then the  result would follow from  this long exact
sequence
since
$\pi_1(\cM^*)=0$.  This is true, in fact the results of \cite{daska},
though stated for $U(2)$ instead of $SU(2)$, are equally valid for
$SU(2)$.

Alternatively, we can show that  $ \pi_1(\cF^*)$ is zero using the
result
of \cite{rade}. Fix a flat irreducible connection $A\in \cF^* $. Let
$\gamma:S^1\to \cF$ be a loop of flat connections. Corollary 2.7 of
\cite{daskalopoulos-uhlenbeck}  shows that
$\cA^s $ is simply connected, and since
$\cF^*\subset \cA^s$ (see \cite{Atiyah-Bott})  
$\gamma$ bounds a disc
$\bar{\gamma}:D^2\to \cA^s $. The main result of \cite{rade}
asserts that
the Yang-Mills flow on $\cA$, restricted to
$\cA^s$, defines a strong deformation retraction of $\cA^s$ to
$\cF^*$. Thus $\gamma$ bounds a disc in $\cF^*$, completing the proof.

\end{proof}

The space of  irreducible $SU(2)$ representations of a free group is
easily proven to be path connected. Thus the space
$\cM^*(X_1)$, homeomorphic to the space of conjugacy classes of
irreducible representations of $\pi_1(X_1)$,   is path connected. There
is
a fibration
$$\cG(X_1)/\pm \Id\to\cF^*(X_1)\to \cM^*(X_1)$$ as in
the
proof of Lemma \ref{lemma6.1}. The gauge group
$\cG(X_1)=\text{Map}(X_1,SU(2))$ is path connected since $X_1$ has the
homotopy type of a 1-complex. Thus $\cF^*(X_1)$ is path connected.
The same argument shows that
    the space of irreducible perturbed flat connections on
$X_2$ is
    also path connected, since   the  perturbation was chosen small
enough
so that $\cM^*_f(X_2)$ is diffeomorphic to the unperturbed
moduli space $\cM^*(X_2)$.

Combining the facts of the preceding paragraph, Lemma \ref{lemma6.1},
and the observation that
since the space of connections on a manifold is contractible, any family
of connections on a submanifold of
$X$ extends to a family of connections on all of $X$,  we obtain the
following useful proposition.

\begin{prop} Given perturbed flat irreducible $SU(2)$ connections $A_0,
A_1$ on the closed 3-manifold $X$, there exists a 2-parameter family
$A_{s,t},\ s,t\in[0,1]$    of irreducible connections on $X$ satisfying:
\begin{enumerate}
\item $A_{0,t}=A_0$  and $A_{1,t}=A_1$ for all
$t\in[0,1]$.
\item $ A_{s,0} $ restricts to a path  of irreducible flat connections
on
$X_1$.
\item $A_{s,1}$ restricts to a path of irreducible perturbed flat
connections on
$X_2$.
\item $A_{s,t}$ restricts to an irreducible  flat connection
$a_{s,t}$ on
$\Sigma$ for all  $s,t$.
\item The restriction of $A_{s,t}$ to a collar of $\Sigma$ is in product
form:
$A_{s,t}|_{\Sigma\times[-1,1]}=\tfrac{\del}{\del u} + a_{s,t}$.
\qed\end{enumerate}
\end{prop}

We now give an argument that $\Sign_T=\Sign_C$.  The strategy is to
identify
$\Sign_T$ with an infinite--dimensional Maslov index of the Cauchy data
spaces, then
to use  homotopy and stretching arguments to symplectically reduce to a
finite--dimensional Maslov index which we identify with $\Sign_C$. We
recommend that
the reader who is unfamiliar with the type of argument we use look at
the ``user's
guide'' section of
\cite{daniel-kirk}.

Let $\La^1_{s,t}$ denote the  Cauchy data space for the restriction of
the
operator
$D_{A_{s,t},f}$ to $X_1$. Similarly define $\La^2_{s,t}$. These vary
continuously by Theorem \ref{thm5.3}.  Note that in contrast to Section 3, the superscripts 
$1$ and $2$ refer to the Cauchy data spaces for the two parts of the Heegard decomposition, 
not to the Sobolev indices. Thus $\La^1_{s,t}$ and $\La^2_{s,t}$ are continuous families of Lagrangian subspaces of $L^2(\Om^*_\Sigma(\ad P))$. Moreover, for each $(s,t)$, $(\La^1_{s,t},\La^2_{s,t})$ is  a Fredholm pair.  Hence the Maslov index $\Mas(\La^1_\alpha,\La^2_\al)\in\zz$ is defined
for any path $\al:[0,1]\to [0,1]\times[0,1]$ (for details see \cite{kirk-lesch}).

    Nicolaescu's theorem  says that the spectral flow of the family
$D_{A_{s,0},f},
\ s\in[0,1]$ equals the (infinite--dimensional) Maslov index
$\Mas(\La^1_{s,0},\La^2_{s,0})$ (see \cite{nico,kirk-lesch}).
Hence
$$\Sign_T([A_0],[A_1])=(-1)^{\SF(D_{A_{s,0}})}
=(-1)^{\Mas(\Lambda^1_{s,0},\Lambda^2_{s,0})}.$$
We will
use
the homotopy properties of the Maslov index to find a simpler,
finite--dimensional expression  for $\Mas(\La^1_{s,0},\La^2_{s,0})$ which
can  be compared to the difference in intersection numbers.

The reduction  will be based on the following observations. Since the
connection $a_{s,t}$ on $\Sigma$ is flat for all $s,t$,  the kernel of
$S_{a_{s,t}}$ is identified via the Hodge theorem with  the cohomology
$H^*(\Sigma;\ad P_{a_{s,t}})$ of the twisted de Rham complex
\begin{equation}\label{eq6.10}  0\to \Omega^0(\Sigma,\ad
P_{a_{s,t}})\xrightarrow{d_{a_{s,t}}}
\Omega^1(\Sigma,\ad P_{a_{s,t}})\xrightarrow{d_{a_{s,t}}}
\Omega^2(\Sigma,\ad P_{a_{s,t}})\to 0.
\end{equation}

Since $a_{s,t}$ is irreducible,
$H^0(\Sigma;\ad P_{a_{s,t}})=0$ (i.e.~using \eqref{eq2.3} for $\Sigma$)
and so
$H^2(\Sigma;\ad P_{a_{s,t}})=0$ by Poincar\'e duality.  Thus the kernel
of
$S_{a_{s,t}}$ is just
$H^1(\Sigma;\ad P_{a_{s,t}})$. But this   is canonically identified
with the tangent space of
$\cM^*(\Sigma)$ at $[a_{s,t}]$, since the image of
$d_{a_{s,t}}:\Omega^0(\Sigma,\ad P_{a_{s,t}})\to
\Omega^1(\Sigma,\ad P_{a_{s,t}})$ is the tangent space of the
$\cG(\Sigma)$
orbit and   $d_{a_{s,t}}:\Omega^1(\Sigma,\ad P_{a_{s,t}})\to
\Omega^2(\Sigma,\ad P_{a_{s,t}})$ is the linearization at $a_{s,t}$ of
the
curvature map
$a\mapsto F(a)$.

    Via this identification, the cup product (or wedge product if one
identifies cohomology with  harmonic forms)
$$H^1(\Sigma;\ad P_{a_{s,t}})\times H^1(\Sigma;\ad P_{a_{s,t}})
\xrightarrow{\cup}H^2(\Sigma;\rr)=\rr$$ defines a symplectic structure
on
$\cM^*(\Sigma)$. In particular, since
$\cM^*(\Sigma)$ is a smooth $6g-6$ dimensional manifold,
$\ker S_{a_{s,t}}=H^1(\Sigma;\ad P_{a_{s,t}})$ is a smoothly varying
family of symplectic subspaces of $L^2(\Omega^*_\Sigma(\ad P))$.  For
convenience we use the notation
$$H^1_{s,t}(\Sigma)=H^1(\Sigma;\ad P_{a_{s,t}}).$$

Let $P^+_{s,t}$ denote the ($L^2$ closure of the) positive eigenspan of
the
tangential operator
$S_{a_{s,t}}$, and let  $P^-_{s,t}$ denote its negative eigenspan.
Thus
we have a smoothly varying (in $s,t$) decomposition
$$L^2(\Omega^*_\Sigma(\ad P))=P^-_{s,t}\oplus H^1_{s,t}(\Sigma)\oplus
P^+_{s,t}.$$

\begin{lemma} \label{lem6.1} For each $s\in [0,1]$, let
$T^1_{s,0}\subset
H^1_{s,0}(\Sigma)$  be the   subspace defined  by
$$T^1_{s,0}=\text{\rm image }H^1(X_1;\ad P_{a_{s,0}})\to H^1(\Sigma;\ad
P_{a_{s,0}}).$$ Then the spaces $T^1_{s,0}$ are Lagrangian  and vary
continuously. Via the identification  $T_{a_{s,0}}\cM^*(\Sigma)\cong
H^1_{s,0}(\Sigma)$ the subspace $T^1_{s,0}$ corresponds to the tangent
space at $a_{s,0}$ of the submanifold $\cM^*(X_1)$ of $\cM^*(\Sigma)$.

Moreover
$$\SF(D_{A_{s,0},f})=\Mas(\La^1_{s,0},\La^2_{s,0})=
\Mas(P^+_{s,0}\oplus T^1_{s,0},\La^2_{s,0}).$$
\end{lemma}

\begin{proof}The statement of Lemma \ref{lem6.1}  refers only to the
parameters $(s,0), s\in [0,1]$; in other words, throughout the proof the
``$t$'' parameter is zero.

The fact  that $T^1_{s,0}$ is Lagrangian is a  standard fact in
geometric
topology whose easy proof (based on Poincar\'e duality) we leave to the
reader.  The identification of
$T^1_{s,0}$ with the tangent space of $\cM^*(X_1)$ is also standard; one
of many proofs is to consider the map of  de Rham complexes
$\Omega^*(X_1;\ad P)\to \Omega^*(\Sigma;\ad P)$ in light of the diagram
\eqref{eq2.3}.
   Since $\cM^*(X_1)$ is a
smooth submanifold, it follows that $T^1_{s,0}$ varies smoothly. Thus we
focus on the statements about the Maslov index.

It follows from the adiabatic limit theorem of \cite{nico} and its
refinement in
\cite{daniel-kirk} that   replacing the collar of $X_1$ by an
increasingly
longer collar $\Sigma\times [-R,0]$ gives a homotopy of the path
$s\mapsto
\La^1_{s,0}$ to
$s\mapsto P^+_{s,0}\oplus T^1_{s,0}$.  This is a consequence of the
fact that
there are no
$L^2$ solutions to $D_{A_{s,0},f}\phi=0$ on  $(\Sigma\times
(-\infty,0])\cup
X_1$, which follows from \cite{APS} and two facts:
\begin{enumerate}
\item  The  restriction map
$H^1(X_1;\ad P_{a_{s,0}})\to H^1(\Sigma;\ad P_{a_{s,0}}) $  is
injective.
\item The space $T^1_{s,0}$  equals  the {\em limiting values of
extended
$L^2$ solutions} in the sense of \cite{APS}. This was first observed by
Yoshida \cite{yoshida} and in the present notation a proof can be found
in
\cite{daniel-kirk} and \cite{kirk-lesch}.
\end{enumerate}

Let $r\mapsto \La^1_{s,0}(r)$ denote this homotopy. To be explicit,
\[
\La^1_{s,0}(r)=\begin{cases}\text{the Cauchy data space for   }
D_{A_{s,0},f}
\text{ on }(\Sigma\times [-\tfrac{1}{1-r},0])\cup X_1&\text{if } r<1,\\
P^+_{s,0}\oplus T^1_{s,0}&\text{if } r=1.
\end{cases}
\]

This is not a homotopy rel endpoints, and so to see why the Maslov
indices
nevertheless agree, it suffices to show that  the dimensions of the
intersections  $\La^1_{0,0}(r)\cap
\La^2_{0,0}$ and
$\La^1_{1,0}(r)\cap
\La^2_{1,0}$ are independent of $r$, in fact these intersections are
zero
for
all $r$.

For $r<1$ these intersections are isomorphic to the kernel of
    the extension of $D_{A_{i,0},f}$ $(i=0,1)$ to $X(r)=(\Sigma\times
[-\tfrac{1}{1-r},0])\cup X_1\cup X_2$, i.e. $X$ with a long collar
inserted near $\Sigma$.  But the dimension of the kernel of
$D_{A_{i,0},f}$ is independent of the choice of  Riemannian metric on
$X$
since $A_{i,0}$ is perturbed flat: it is isomorphic to the cohomology of the
complex
\eqref{eq2.6}, which is zero by the assumption that each perturbed flat
connection is non-degenerate.

For $r=1$, we need to see that $(P^+_{0,0}\oplus T^1_{0,0})\cap
\La^2_{0,0}$ and
$(P^+_{1,0}\oplus T^1_{1,0})\cap \La^2_{1,0}$ are both zero. The
quotient
$$\frac{(P^+_{0,0}\oplus \ker S_{a_{0,0}}) \cap \La^2_{0,0}}{{
{P^+_{0,0} \cap
\La^2_{0,0}}}}\subset \ker S_{a_{0,0}} $$ is just $T^2_{0,0}$
(\cite{APS,kirk-lesch})
and the fact that each perturbed flat connection is non-degenerate
implies
that $T^1_{0,0}\cap T^2_{0,0}=0$. Moreover,  $ P^+_{0,0} \cap \La^2_{0,0}=0$ since 
any section in this intersection is the restriction to the boundary of a solution to $D_{A_{0,0},f}\phi=0$ which extends to an $L^2$ solution on 
$(\Sigma\times
(-\infty,0])\cup
X_1$.   The only such $\phi$ is zero, as we argued above. 
Hence
$(P^+_{0,0}\oplus T^1_{0,0})\cap \La^2_{0,0}=0$. The other endpoint is
handled
identically.

Since  the Maslov index is unchanged when the paths of Lagrangians are
homotoped by a
homotopy that preserves the dimension of the intersections at the
endpoints
(\cite{kirk-lesch}),
$$\Mas( \La^1_{s,0},\La^2_{s,0})=\Mas(P^+_{s,0}\oplus
T^1_{s,0},\La^2_{s,0}),$$ as
desired.
\end{proof}

The spaces $H^1_{s,t}=\ker S_{s,t}$ form a symplectic vector bundle over
the square
$[0,1]\times[0,1]$. In fact, this is just the pullback of the tangent
bundle  $T_*\cM^*(\Sigma)$ via the map $[0,1]\times [0,1]\to
\cM^*(\Sigma)$
given by $(s,t)\mapsto [a_{s,t}]$.

    We have a Lagrangian subbundle
$T^1_{s,0}$ defined over the top edge of this square, and also the 2
sides
since
$A_{0,t}$ is independent of $t$ and $A_{1,t}$ is independent of $t$.
Extend
this arbitrarily (e.g.~ pull back via a retraction  of the square to
its
three sides) to a Lagrangian subbundle
$T^1_{s,t}$ over the entire square. Then the 2-parameter family
$P^+_{s,t}\oplus T^1_{s,t}$  provides a homotopy which shows that
\begin{equation}\label{eq6.1}
\Mas(P^+_{s,0}\oplus T^1_{s,0},\La^2_{s,0})=\Mas(P^+_{s,1}\oplus
T^1_{s,1},\La^2_{s,1}).
\end{equation} Together with Lemma \ref{lem6.1}, \eqref{eq6.1} implies
that
\begin{equation}\label{eq6.2} SF(D_{A_{s,0},f})= \Mas(P^+_{s,1}\oplus
T^1_{s,1},\La^2_{s,1}).
\end{equation} An argument very similar to the proof of Lemma
\ref{lem6.1}, this time taking the adiabatic limit on $X_2$,  shows that
\begin{equation}\label{eq6.3}
     \Mas(P^+_{s,1}\oplus T^1_{s,1},\La^2_{s,1})=\Mas(P^+_{s,1}\oplus
T^1_{s,1},P^-_{s,1}\oplus T^2_{s,1}),
\end{equation} using the fact that the path of connections $A_{s,1}$ is
perturbed flat on $X_2$.  Moreover, $T^2_{s,1}\subset H^1_{s,1}(\Si)$ is
isomorphic to image of the differential of the inclusion
$\cM_f^*(X_2)\to \cM^*(\Sigma)$ at $A_{s,1}$, i.e.~
$T^2_{s,1}=\text{image}\big( T_{A_{s,1}}(\cM_f^*(X_2))\to
T_{a_{s,1}}(\cM^*(\Sigma))\big)$.
       Combining \eqref{eq6.2} and \eqref{eq6.3} and the facts that  the
Maslov index is additive with respect to direct sum of symplectic spaces
and that $P^+_{s,1}$ and $P^-_{s,1}$ are always orthogonal, one arrives
at
the conclusion
\begin{equation}
\label{eq6.5}
\SF(D_{A_{s,0},f})=\Mas( T^1_{s,1},  T^2_{s,1}).
\end{equation}

The spectral flow $\SF(D_{A_{s,0},f})$ depends only on $f$ and the
endpoints
$A_0$ and
$A_1$ of $A_t$; see  Lemma
\ref{independence}. Moreover its reduction modulo 2 (in fact modulo 8)
only
depends on the gauge equivalence classes $[A_0]$ and $[A_1]$.
Equation \eqref{eq6.5} has identified $\SF(D_{A_{s,0},f})$
  with a finite--dimensional Maslov index  corresponding to
the tangent spaces of Lagrangian submanifolds of  a symplectic
manifold. The
following   lemma  identifies the mod 2 reduction of this Maslov index
with
the difference in intersection numbers of  $\cM^*(X_1)$ and
$\cM^*(X_2)$ at
$[A_0]$ and $[A_1]$.

    \begin{lemma}\label{lem6.2} Let $A$ and $B$ be  smooth  oriented
Lagrangian submanifolds of  a (finite--dimensional) symplectic manifold
$W$.   Let $p$  and $q$ be two transverse intersection points of
$A$ and $B$. Let $\gamma_A$ be a path in $A$ from $p$ to $q$ and
$\gamma_B$ a path in
$B$ from $p$ to  $q$.  Suppose that  $F(s,t)$ is a homotopy rel
endpoints
in $W$ from
$\gamma_A$ to  $\gamma_B$.  Let $T^A_{s,0}$ be the Lagrangian subbundle
of
the restriction of
$T_*W$ to $\gamma_A$ consisting of tangent vectors to $A$, and similarly
let
$T^B_{s,1}$ be the Lagrangian subbundle of the restriction of
$T_*W$ to $\gamma_B$ consisting of tangent vectors to $B$.

Then $T^A$ extends over the homotopy, and if $T^A_{s,t}$ is any
extension
(with
$T^A_{0,s}$ and $T^A_{1,t}$ independent of $t$),  then
$\Mas(T^A_{s,1},T^B_{s,1})$  is even if and  only if  intersection
number
of $A$  and $B$ at $p$ and $q$ coincide.
\end{lemma}\qed

The straightforward proof can easily be constructed by musing on the
following picture, and is left to the reader.
\medskip
\begin{center}
\psfrag{TAs,t}{$T^A_{s,t}$}
\psfrag{(s,t)}{$(s,t)$}
\psfrag{A}{$A$}
\psfrag{B}{$B$}
\psfrag{p}{$p$}
\psfrag{q}{$q$}
\includegraphics{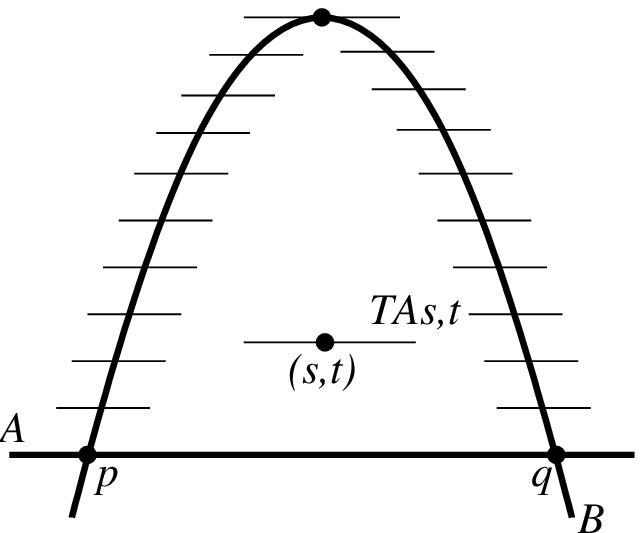}
\end{center}
\medskip
Lemma \ref{lem6.2} and Equation \eqref{eq6.5}  then imply Taubes'
theorem.
\begin{thm}[Taubes]$\Sign_C=\Sign_T$, and hence Casson's   and
Taubes' invariants are equal up to an overall sign (which may depend
on the
homology 3-sphere).\end{thm}\qed

   Taubes proves more, namely he gives a construction  for specifying an
overall sign
so that  in fact  his invariant equals Casson's invariant for all
homology 3-spheres.

\vfill\eject

\end{document}